\documentclass[final,faculty=fw,department=wis,phddegree=mat]{adsphd}

\title{Toric degenerations and Newton-Okounkov bodies}
\author{Cas}{Proost}
\promotor{Prof.~Dr.~F.~Mohammadi}{}
\supervisor{F.~Zaffalon}{}

\researchgroup{Algebra}
\email{cas.proost@student.kuleuven.be} 

\address{Celestijnenlaan 200B}

\date{Academic year 2023-2024}
\copyyear{2024}

\usepackage{textcomp} 
\usepackage{pifont}   
\usepackage{booktabs}
\usepackage[english]{babel}
\usepackage{amsmath, amssymb,amsthm, tikz,  caption}
\usepackage{cite}
\usepackage{graphicx}
\graphicspath{ {./images/} }

\theoremstyle{plain}
\newtheorem{theorem}{Theorem}
\newtheorem{prop}{Proposition}
\newtheorem{lemma}{Lemma}
\theoremstyle{definition}
\newtheorem{defn}{Definition}
\newtheorem{ex}{Example}
\newtheorem{remark}{Remark}
\newtheorem{notation}{Notation}

\usepackage{tikz-cd}
\usepackage{bbm}
\usepackage{float}
\usepackage[noend,ruled,vlined]{algorithm2e}

\usepackage{cleveref}
\crefname{ex}{Example}{Examples}
\crefname{algorithm}{Algorithm}{Algorithms}
\crefname{figure}{Figure}{Figures}
\newcommand{\N}{\mathbb{N}}
\newcommand{\Z}{\mathbb{Z}}
\newcommand{\Q}{\mathbb{Q}}
\newcommand{\R}{\mathbb{R}}
\newcommand{\C}{\mathbb{C}}
\newcommand{\Proj}{\mathbb{P}}
\newcommand{\initial}{\text{in}}
\DeclareMathOperator{\trop}{trop}
\DeclareMathOperator{\Gr}{Gr}
\DeclareMathOperator{\relint}{relint}
\newcommand\scalemath[2]{\scalebox{#1}{\mbox{\ensuremath{\displaystyle #2}}}}

\begin{document}


\makefrontcoverXII

\maketitle

\frontmatter 

\includepreface{preface}
\includeabstract{abstract}



\tableofcontents


\mainmatter 


%
  \graphicspath{{\chaptersdir/introduction/image/}}%
\chapter{Introduction}\label{ch:introduction}




Algebraic varieties are families of solutions to systems of polynomial equations. The vast field of algebraic geometry studies these varieties and, in doing so, has led to many innovations in many areas of science, as polynomial equations are omnipresent. In general however, studying algebraic varieties is difficult.

Some varieties are easier to study than others. One such class are toric varieties, studied in toric geometry. Toric varieties (see \cite{cox2011toric}, \cite{fulton1993introduction}) are a well-understood class of algebraic varieties that are defined by binomial equations. They are equipped with an expansive array of combinatorial tools that allow for explicit calculation of many of their invariants such as Hilbert polynomials or cohomology groups.

A toric degeneration is a mathematical construction that allows for the use of the combinatorial tools from toric geometry to study general algebraic varieties. They can be thought of as deformations of a variety into a toric variety. This allows for the use of those combinatorial tools, on the (not necessarily toric) variety. The innovative and rapidly evolving theory of toric degenerations has proven useful in many any areas of pure mathematics \cite{bossinger2023survey}, such as in algebraic geometry \cite{COLLART1997443}, symplectic geometry \cite{nishinou2010potential}, representation theory \cite{representation}$\ldots$, as well as in applied mathematics, for example in algebraic statistics \cite{Quasisymmetry}, in theoretical physics (mirror symmetry \cite{batyrev1998mirror}, scattering amplitudes in particle physics \cite{scattering} ...), in mechanical engineering \cite{robotics} and even in chemistry \cite{millan2011chemical}, \cite{faulstich2024algebraic}.

However, obtaining toric degenerations of a specific variety is not a trivial problem. In this thesis, a connection is explored between two different areas of mathematics that can be used to obtain toric degenerations: tropical geometry and the theory of Newton-Okounkov bodies. 

One way to construct toric degenerations of a variety is through the process of tropicalization. The tropicalization of a variety is a combinatorial silhouette of the variety. Tropicalizations, or tropical varieties, are the main objects studied in the field of tropical geometry \cite{IntroToTropGeom}. They are polyhedral fans, which are simplicial complexes of linear half-spaces called polyhedral cones. In general they are constructed from tropical semi-rings, which replace the usual addition operation with taking the minimum of the terms, and replace the usual multiplication operation with addition. 
Tropical geometry is a much broader field than toric degenerations and, as such, has an even wider variety of applications \cite{ParallelComputing}, for example in the additional fields of optimization theory \cite{optimization}, economics \cite{tran2019product} and biology \cite{Yoshida2017TropicalPC}. 

Newton-Okounkov bodies are the second way of constructing toric degenerations that we consider. They can be seen as generalisations of the Newton polytope of a polynomial, introduced by Russian mathematician Andrei Okounkov. Their connections to toric degenerations were first introduced in \cite{Anderson}, and expanded upon in \cite{KhovanskiiBases}. 
They are convex bodies that arise when considering valuations on algebras, in particular on coordinate rings, and in a way are less polyhedral analogues for general varieties of the combinatorial objects related to toric varieties. For example, the degree of a projective variety is equal to the (normalized) volume of one of its Newton-Okounkov bodies \cite[Theorem 2.23-2.24]{KhovanskiiBases} (assuming the coordinate ring is contained in a finitely generated algebra and the valuation has one-dimensional leaves).  

One of the families of varieties that is often studied through toric degenerations is that of Grassmannians. The Grassmannian $\Gr(k,n)$ is the set of $k$-dimensional linear subspaces of $n$-space, and effectively generalizes projective space. They have intricate connections to 
combinatorics, which make them very interesting to look at through the use of toric degenerations. In particular, the Grassmannians of planes, i.e.  with $k=2$, are well understood and are suitable for beautiful explicit examples such as in section $5$ of \cite{WallCrossing}. Higher dimensional Grassmannians however, still pose significant challenges, and are the subject of much research in tropical geometry \cite{ParallelComputing} and toric degenerations  \cite{ToricdegensofGrassmannians}.

These challenges partially lie in the fact that not any cone in the tropicalization of a variety or any valuation on the coordinate ring of that variety yields a toric degeneration.
Only specific cones and valuations give rise to toric degenerations and, in general, finding these is not an easy task. 
Indeed, for high dimensional varieties, the amount of cones in a tropicalization can grow exponentially, and the fraction of prime cones, which are the cones that yield toric degenerations, seemingly decreases. 
Attempts have been made to utilise powerful computing techniques in order to compute tropicalizations, for example in \cite{ParallelComputing} parallel computing is used to compute the tropical Grassmannian $\trop \Gr(3,8)$. This still required the intelligent use of various symmetries of $\Gr(3,8)$, and higher dimensional tropical Grassmannian have thus far resisted any attempts at computation.

In \cite{KhovanskiiBases}, the authors prove the existence of a link between the prime cones of a tropical variety and Newton-Okounkov polytopes of the variety. Building on this result, Harada and Escobar wrote the wall-crossing paper \cite{WallCrossing} 
which serves as the foundation for this thesis. 
In the above-mentioned paper, the authors used the previously shown link between prime cones and Newton-Okounkov polytopes to prove that traveling between adjacent prime cones of a tropical variety corresponds to so-called mutations between their associated Newton-Okounkov polytopes. Mutations of polytopes are maps between polytopes that arise in many areas of combinatorics and were introduced in context of the classification of Fano varieties  \cite{Akhtar_2012}, \cite{coates2022mirror}. 
In some situations, these mutations of polytopes can be described in terms of associated combinatorics, such as the rotations of edges in phylogenetic trees for $\Gr(2,n)$.

The wall-crossing theorem effectively allows one to traverse a tropical variety from one prime cone to the next, crossing the walls between the cones, if one knows what form the mutations take. Not only does this mean new toric degenerations can be constructed, but this also allows for computation of points of tropical varieties. However, this only allows wall-crossing between adjacent prime cones. In particular, prime cones that are isolated and not adjacent to other prime cones would be unreachable.

This thesis aims to work towards a generalization of the wall-crossing theorem that allows for a wall-crossing through non-prime cones. The ultimate goal is to compute and study tropicalizations of higher dimensional Grassmannians, by traversing them. For each Grassmannian, there is a canonical Newton-Okounkov polytope, called the Gelfand-Zetlin polytope, which can serve as a starting point for traversing the tropical variety \cite{KOGAN20051}. 

The structure of the thesis is as follows: \Cref{ch : preliminaries} introduces basic preliminaries, from the Grassmannians, to Gr\"obner theory and tropical geometry. \Cref{ch:Tdegenerations} discusses toric degenerations, how they can be obtained from tropicalizations and from Newton-Okounkov bodies, and the connection between the two, including the wall-crossing theorem for prime cones. \Cref{ch:Wall-crossing for non-prime cones} contains our results for non-prime cones, where we propose a way to work around the non-prime cones and apply our approach to the Grassmannian $\Gr(3,6)$ using the computer algebra system Macaulay2 \cite{M2}, showing that wall-crossing through non-prime cones is possible.

\cleardoublepage

%

%
  \graphicspath{{\chaptersdir/Preliminaries/image/}}%
  \chapter{Preliminaries}\label{ch : preliminaries}

This chapter aims to provide the reader with a basic understanding of the underlying theory of the thesis. It also serves to introduce notation and provide some insight into tropical Grassmannians and challenges associated with computing them. The sections on tropical geometry and the Grassmannians are mainly based on \cite{IntroToTropGeom}, while the section on Gr\"obner theory follows \cite{GrobnerBasisAndConvexPolytopes} and the section on polyhedral geometry follows chapter 7 \& 8 of \cite{Schrijver1986TheoryOL}.

\section{The Grassmannians}
This section gives a basic introduction to the family of varieties which serves as the motivating examples for this thesis: the Grassmannians. An overview on Grassmannians can be found, for example, in \cite{IntroToTropGeom}, which also introduces them in the context of tropical geometry.

Grassmannians are a family of varieties that can be thought of as generalizations of projective spaces. Recall that a projective space parametrizes all one dimensional vector subspaces of some ambient space. For example: points in the projective space $\R\Proj^2$ correspond to lines through the origin in $\R^{3}$. Grassmannians parametrize all linear subspaces of a fixed dimension. 

\begin{defn}[Grassmannians]
    The Grassmannian $\Gr(k,n)$ is the set of $k$-dimensional vector sub-spaces of an $n$-dimensional vector space.
\end{defn}

\begin{remark}
    Grassmannians can be defined over any field, but in this thesis they will be considered mostly over $\C$. The main exception will be when talking about higher dimensional tropical Grassmannians and their dependence on the characteristic of the field. 
\end{remark}

It is uncommon to see an author use the notation $\Gr(1,n+1)$, as this is the $n$-dimensional projective space and it is therefore usually denoted $\Proj^n$. Another convention is that when talking about Grassmannians, we consider $k\le n/2$, since there is a duality between $k$-dimensional subspaces of $n$-space and $n-k$ dimensional subspaces, and hence we can identify $\Gr(k,n) \cong \Gr(n-k,n)$.

To represent a point in projective space, the notion of homogeneous coordinates is used. A point in projective space is not represented by a single coordinate, but rather by an equivalence class of coordinates, where two coordinate vectors are equivalent if one is a (nonzero) multiple of the other. Indeed any nonzero point in $n$-space determines the line through the origin and that point, but any other point on that line would yield the same line. For Grassmannians, one needs $k$ linearly independent points of $n$-space to represent a point of $\Gr(k,n)$. This gives a $k\times n$ matrix.
It is clear that a point of $\Gr(k,n)$ does not have a unique representing matrix but, unlike the case of projective spaces, it is not immediately clear when two matrices represent the same point. To answer this question we need to introduce a map called Pl\"ucker embedding. 

The Pl\"ucker embedding is a map embedding $\Gr(k,n)$ into a $\left({{n}\choose{k}}-1\right)$-dimensional projective space, by sending such a matrix representative of a point in $\Gr(k,n)$ to the vector of its $k\times k$ minors.
These minors are called the Pl\"ucker coordinates or variables. We introduce some notation in order to formally define this embedding.

\begin{notation}
    We will use the notations $[n]:= \{1, \ldots, n\}$, and \[ {{[n]}\choose{k}} := \{I \vert I \subset [n], |I| = k\}.\]
\end{notation}

For a $k\times n$ matrix $M$, denote by $M_I$ the sub-matrix of $M$ obtained by 
selecting the columns of $M$ indexed by $I \in {{[n]}\choose{k}}$ and arranging them in  order of increasing index.

\begin{defn}[Pl\"ucker variables]
    We label the coordinates on $\Proj^{{n\choose k} -1}$ with the elements of ${[n]\choose k}$. The coordinate indexed by an $I\in {[n]\choose k}$ is called the Pl\"ucker variable or Pl\"ucker coordinate $p_I$. 
\end{defn}
\begin{remark}\label{plucker permutation}
    Occasionally, we will write a Pl\"ucker coordinate whose indices are not in increasing order: suppose $J$ is a permutation of $I\in {[n]\choose k}$. Then we will use the convention $p_J = \mathrm{sgn}(J)p_I$. For example $p_{21} = -p_{12}$.
\end{remark}

\begin{defn}[Pl\"ucker embedding]
    Let $M$ be a matrix representing a point $x$ of $\Gr(k,n)$. The Pl\"ucker embedding is the map 
    \[\Gr(k,n)\to \Proj^{{n\choose k}-1}\]
    which maps $x$ to the point in $\Proj^{{n\choose k}-1}$ with $p_I = \det(M_I)$ for all $I\in {[n]\choose k}$.
\end{defn}

\begin{remark}
    This map also makes sense with respect to the convention in \cref{plucker permutation}: if $I \in {[n]\choose k}$  and $J$ is a permutation of $I$, then $|M_J| = \mathrm{sgn}(J)|M_I|$ where $M_J$ is the matrix with the columns put in the order described by $J$.
\end{remark} 

The Pl\"ucker embedding has an important property: for any given $k\times n$ matrix, the coordinates of its image under the Pl\"ucker embedding
satisfy some algebraic relations, called the Pl\"ucker relations, that follow from linear algebra. Under the Pl\"ucker embedding, $\Gr(k,n)$ is exactly the projective variety of the ideal generated by all Pl\"ucker relations. Herein lies the power of the Pl\"ucker embedding: it allows the Grassmannians to be viewed as a projective variety instead of an abstract set of linear subspaces 
. 
Since the Grassmannians are most interesting as a variety, many authors define them by use of this Pl\"ucker embedding instead of as the set of k-dimensional subspaces of some vector space. 

We give an explicit form of the Pl\"ucker relations. Take $I,J \in {{[n]}\choose{k}}$ disjoint, and let $i_a \in I$ be some coordinate of $I$, then denote by $I(i_a \to j)$ the set obtained by replacing $i_a$ in $I$ by $j$, and analogously, denote by $J(j\to i_a)$ the set obtained by replacing $j$ in $J$ by $i_a$. A Pl\"ucker relation is of the form:

\[p_Ip_J = \sum\limits_{j\in J} p_{I(i_a \to j)}p_{J(j\to i_a)}.\]

To obtain all Pl\"ucker relations, it is sufficient to only vary $I$ and $J$, and not $i_a$. For example one can pick $i_a$ to always be the final coordinate of $I$. Nevertheless writing all the Pl\"ucker relations down is tedious to do by hand, especially for higher $k$ and $n$.

\begin{defn}[Pl\"ucker ideal]
    The Pl\"ucker ideal $I_{k,n}$ is the ideal, in the ring $\C[p_I|I\in{{[n]}\choose{k}}]$ of polynomials in the Pl\"ucker variables, generated by all the Pl\"ucker relations of $\Gr(k,n)$.
\end{defn}

\begin{ex}
    The smallest Grassmannian that is not a projective space is $\Gr(2,4)$. Recall that points in this Grassmannian correspond to planes in four dimensional space. Each such plane can be represented by two linearly independent vectors with four coordinates, hence by a ($2\times 4$) matrix.

    The Pl\"ucker variables are indexed by the elements of ${[4]\choose 2}$ and therefore the complete set of Pl\"ucker coordinates is: 
    \[\{p_{12}, p_{13}, p_{14}, p_{23}, p_{24}, p_{34}\}.\]
    Let us work over $\R$ for simplicity, and consider the plane in $\R^4$ generated by vectors $(4,3,2,1)$ and $(1,2,3,4)$. This plane is an element of the real Grassmannian $\Gr(2,4)$, and we can represent it 
    as the span of the rows of the matrix:
    \[M = \begin{bmatrix}
        4 & 3 & 2 & 1\\
        1 & 2 & 3 & 4
    \end{bmatrix}\]
    Its Pl\"ucker coordinate $p_{12}$ is given by the determinant of the submatrix of $M$ consisting of the first column followed by the second column:
    \[\begin{vmatrix}
        4 & 3\\
        1 & 2    
    \end{vmatrix}  = 5.\]
    Similarly, the other Pl\"ucker coordinates can be obtained to find the vector $(5, 10, 15, 5, 10, 5)$, where the order is as given before. Note that the Pl\"ucker embedding is into projective space, so these coordinates are homogeneous, and the same element of $\Gr(2,4)$ is represented by the Pl\"ucker coordinates $(1,2,3,1,2,1)$.\\
    Let us now find the Pl\"ucker relations for $\Gr(2,4)$. Let us start with $I = \{1,2\} $ and $J = \{3,4\}$. Letting $i_a$ be the last coordinate of $I$ as before we get the relation
    \[p_{12}p_{34} = p_{13}p_{24} + p_{14}p_{32}\]
    which using $p_{23} = -p_{32}$ we can rewrite as
    \[p_{12}p_{34} - p_{13}p_{24} + p_{14}p_{23} = 0.\]
    This is in fact the only Pl\"ucker relation of $\Gr(2,4)$, since taking $I,J$ to be $\{1,3\}, \{2,4\}$ or $\{1,4\}, \{2,3\}$ yields the same relation. Therefore the image of the Grassmannian $\Gr(2,4)$ under the Pl\"ucker embedding is the vanishing locus of the ideal $I_{2,4} = \langle p_{12}p_{34} - p_{13}p_{24} + p_{14}p_{23}\rangle$. \\
    Finally let us check that our Pl\"ucker relation does in fact vanish on our aforementioned element of $\Gr(2,4)$, by filling in the coordinates $(1,2,3,1,2,1)$. We get
    \[1\cdot 1 - 2 \cdot 2 + 3 \cdot 1 = 0\] as expected.
\end{ex}

\section{Polyhedral and tropical geometry}

In this section we introduce the main combinatorial objects of interest to this thesis: tropical varieties. In our context, they can be thought of as ``combinatorial shadows'' of algebraic varieties. Though the theory of tropical geometry is a rich and independent field of combinatorics, we will introduce tropical varieties from the point of view of Gr\"obner theory. 

\subsection{Polyhedral geometry}

We start by giving a concise introduction to the relevant concepts from polyhedral geometry. This allows us to introduce the combinatorial objects of which tropical varieties are examples: polyhedral fans. For a more complete view  of polyhedral geometry the reader is referred to chapter 7 \& 8 of \cite{Schrijver1986TheoryOL}.

\subsubsection{Polytopes and polyhedral cones}

Polyhedral geometry is concerned with polyhedra, which are intersections of a finite amount of ``half-spaces''. Intuitively, one takes a finite amount of hyperplanes and, for each hyperplane, removes all points on one side of the hyperplane. What is left is a polyhedron. 

\begin{ex}
    The simplest example of a polyhedron is a two dimensional polygon, such as a square or a triangle. A circle is not a polyhedron, as it can not be written as a \textit{finite} intersection of half spaces. 
\end{ex}

In the following we will work over $\R$. A hyperplane in $n$-space is given by a linear equation in $n$ variables say 
\[a_1x_1+ \cdots + a_n x_n = b\]
which can be written as $a\cdot x = b $ where $a:= (a_1,\ldots, a_n)$ and $x^\top := (x_1, \ldots, x_n)$. To obtain a half space lying on one side of this hyperplane we replace the equality with an inequality. A finite intersection of half spaces then is given by a finite system of such inequalities, which can be written using matrices as \[Ax \ge b\] where the inequality is term wise between the vectors. 

\begin{defn}[Polytope]
    A polytope is a bounded polyhedron. They can be written as the convex hull of a finite set of points $p_1, \ldots ,p_k$:
    \[\mathrm{conv}(\{p_1,\ldots,p_k\}) := \left\{\lambda_1p_1+\cdots+\lambda_kp_k 
    \vert \; \lambda_i\ge0, \sum_i\lambda_i = 1\right\}\]
\end{defn}
The second part of this definition is a property of polytopes that is frequently used as a definition. 
\begin{ex}
    Points, line segments and polygons are all examples of polytopes. So are three dimensional polyhedra like a cube, tetrahedron or dodecahedron.
\end{ex}

We will not be interested in only bounded polyhedra, but also in a certain class of unbounded polyhedra: polyhedral cones. 
\begin{defn}[Polyhedral cone]
    A polyhedral cone is the intersection of a finite amount of half spaces whose defining hyperplanes contain the origin. Their equations have the form \[Ax\ge0.\]
    A polyhedral cone can be written as the positive hull of a finite set of points $p_1, \ldots p_k$:
    \[\mathrm{cone}(\{p_1,\ldots,p_k\}) := \left\{\lambda_1 p_1+\cdots+\lambda_k p_k 
    \vert\; \lambda_i\ge0\right\}\]
\end{defn}

\begin{remark}
    The origin is a trivial cone and is the only cone that is bounded (assuming we are working over a field of nonzero characteristic).
\end{remark}
\begin{remark}
    The way we defined cones here is sometimes called a closed polyhedral cone by other authors. This is as opposed to a open polyhedral cone which is where the inequalities in the definition would be strict. We will always mean closed cone unless stated otherwise.
\end{remark}

To study all polyhedra it is sufficient to study all polytopes and cones, even though not every polyhedron is either a polytope or a polyhedral cone. This is because of the decomposition theorem, which states that every polyhedron can be decomposed as the Minkowski sum (a polyhedral analogue to the direct sum of vector spaces) of a polytope and a polyhedral cone. The precise statement of this theorem is not relevant for our purposes, but the interested reader can refer to \cite[Theorem 8.5]{Schrijver1986TheoryOL}.

\begin{defn}[Face of a polyhedron]
    A face of a polyhedron $P\subset \R^n$ is a subset of the form:
    \[F = \text{face}_w(P):= \{p \in P \vert p\cdot w \ge q\cdot w, \forall q \in P\}\]
    with $w\in \R^n$ some vector.
\end{defn}

\begin{ex}
    Any polyhedron $P$ is a face of itself, since $\text{face}_0(P) = P$. The faces of a polygon are exactly the polygon itself, its  
    edges and its vertices.
\end{ex}

We state the following proposition which follows from the definition of a face. 

\begin{prop}
    A face of a polyhedron is a non-empty polyhedron and consequently a face of a polytope, resp. a cone, is a polytope, resp. a cone.
\end{prop}
\subsubsection{Polyhedral fans}

The class of combinatorial object that concerns us most is that of polyhedral fans. These are finite sets of cones that intersect nicely, that is, the set of cones forms a simplicial complex.
\begin{defn}
    A polyhedral fan $\Sigma$ is a finite set of polyhedral cones such that the following properties hold:
    \begin{itemize}
        \item For any $C_1,C_2 \in \Sigma$, $C_1\cap C_2$ is a face of both $C_1$ and $C_2$
        \item  For any $C_1,C_2 \in \Sigma$, $C_1\cap C_2 \in \Sigma$
    \end{itemize}
\end{defn}
 The support of a fan is the union of all its cones, and when the support of a fan is equal to the ambient space, the fan is called complete.
 
\begin{ex}\label{example:fan}
     
     Consider the cones $C_1 = \mathrm{cone}(e_1,e_2),C_2 = \mathrm{cone}(e_1,-e_1-e_2)$ and $C_3 = \mathrm{cone}(e_2,-e_1-e_2)$ in $\R^2$ where $e_1 = (1,0)$ and $e_2=(0,1)$. The set of these cones and every intersection of these cones forms a complete fan.
     
     \begin{figure}[H]
    \begin{center}
        \begin{tikzpicture}[x=0.75pt,y=0.75pt,yscale=-0.5,xscale=0.5]

\draw  [draw opacity=0][fill={rgb, 255:red, 74; green, 144; blue, 226 }  ,fill opacity=0.55 ] (50,400) -- (250,0) -- (250,200) -- cycle ;
\draw  [draw opacity=0][fill={rgb, 255:red, 74; green, 144; blue, 226 }  ,fill opacity=0.55 ] (250,200) -- (450,200) -- (50,400) -- cycle ;
\draw  [draw opacity=0][fill={rgb, 255:red, 74; green, 144; blue, 226 }  ,fill opacity=0.55 ] (250,0) -- (450,200) -- (250,200) -- cycle ;
\draw    (250,400) -- (250,0) ;
\draw    (50,200) -- (450,200) ;
\draw [color={rgb, 255:red, 74; green, 144; blue, 226 }  ,draw opacity=1 ][line width=2.25]    (50,400) -- (250,200) ;
\draw [color={rgb, 255:red, 74; green, 144; blue, 226 }  ,draw opacity=1 ][line width=2.25]    (250,200) -- (450,200) ;
\draw [color={rgb, 255:red, 74; green, 144; blue, 226 }  ,draw opacity=1 ][line width=2.25]    (250,200) -- (250,0) ;

\draw (301,130) node [anchor=north west][inner sep=0.75pt]   [align=left] {$\displaystyle C_{1}$};
\draw (261,220) node [anchor=north west][inner sep=0.75pt]   [align=left] {$\displaystyle C_{2}$};
\draw (201,152) node [anchor=north west][inner sep=0.75pt]   [align=left] {$\displaystyle C_{3}$};

\end{tikzpicture}
\caption{\cref{example:fan}}
    \end{center}
\end{figure}
This is easily seen by computing the intersections
     \[C_1\cap C_2 = \mathrm{cone}(e_1), \hspace{2pt} C_1\cap C_3 = \mathrm{cone}(e_2), \hspace{2pt} C_2\cap C_3 = \mathrm{cone}(-e_1-e_2)\] and \[C_1\cap C_2 \cap C_3 = \{(0,0)\}\] and checking that they satisfy the required properties.
     \tikzset{every picture/.style={line width=0.75pt}} 
\end{ex}

\subsection{Tropical geometry}
In this section we give a very limited introduction to tropical geometry, from the point of Gr\"obner theory. The standard reference for tropical geometry is \cite{IntroToTropGeom}, which most of this section is based on. The first part on Gr\"obner theory mostly follows \cite{GrobnerBasisAndConvexPolytopes}.

In tropical geometry, a variety $V(I)$ is ``tropicalized'' to obtain a polyhedral fan $\trop(I)$ called the tropical variety or tropicalization of $I$ (or of $V(I)$).
This fan structure is inherited from the Gr\"obner fan of the ideal, which is a fan encoding all so-called initial ideals. We therefore start with some basic Gr\"obner theory.

\subsubsection{Some Gr\"obner theory}\label{Gr\"obnertheory}

Let $R:= \C[x_1, \ldots, x_n]$ be a polynomial ring
\begin{defn}[Monomial ordering]
    A monomial ordering $\preceq$ on $R$ is a total order on the set of all monomials in $R$ such that the following holds:
for any two monomials $x^{\alpha_1}, x^{\alpha_2} \in R$, if $x^{\alpha_1} \preceq x^{\alpha_2}$, then for any $x^\beta \in R$ we have $x^{\beta}x^{\alpha_1} \preceq x^{\beta}x^{\alpha_2}$
\end{defn}

\begin{ex}
    The lexicographic monomial order on $\C[x_1, \ldots, x_n]$ is defined as follows: $ x^w\preceq_{\text{lex}} x^v$ for $v,w \in \R^n$ if the leftmost non-zero entry of $w-v$ is positive. The reverse lexicographic order $\preceq_\text{revlex}$ can be defined similarly by looking at the rightmost nonzero entry of $w-v$.
\end{ex}

Monomial orderings allow for a measure of importance of the terms of a polynomial, by considering the minimal term with respect to a certain order. This allows one to reduce a complex polynomial with many terms to a single term, called the initial form, or an ideal of complex polynomials to a monomial ideal, called the initial ideal. Gr\"obner theory concerns the properties of these initial ideals and their relations with their original ideals.

Let $\preceq$ be a monomial ordering on a polynomial ring $R = \C[x_1,\ldots,x_n]$ and $f = \sum_{\alpha \in \N^n} c_\alpha x^\alpha \in R$ be a polynomial, where $c_\alpha \in \C$ and $x^\alpha = x_1^{\alpha_1}\cdots x_n^{\alpha_n}$.

\begin{defn}[Initial form]
The initial form of $f$ with respect to $\preceq$ is defined as follows: 
\[\initial_\preceq(f) := \min_{\alpha:c_\alpha \ne0}(c_\alpha x^\alpha)\] 
where the minimum is taken with respect to $\preceq$.
\end{defn}

Note that this minimum exists since $\preceq$ is a total order and is always attained as a polynomial has only finitely many terms. 
Now let $I \subset R$ be an ideal.
\begin{defn}[Initial ideal]
    The initial ideal of $I$ with respect to $\preceq$ is defined as follows: 
    \[\initial_\preceq(I) := \langle \initial_\preceq(f)|f \in I \rangle.\] 
\end{defn}

For a monomial order $\preceq$, the set of monomials not in $\initial_\preceq(I)$ is called the set of standard monomials $\mathcal{S}(\preceq,I)$. A well-known result of Gr\"obner theory (see \cite[Proposition 1.1]{GrobnerBasisAndConvexPolytopes}) states that the projection of the standard monomials of $\initial_\preceq(I)$ onto $\C[x_1, \ldots, x_n]/I$ forms a vector space basis of this quotient ring.

\begin{remark}
    Note that if $I = \langle f_1, \ldots , f_n \rangle$, it is in general not true that $\initial_\preceq(I) = \langle \initial_\preceq(f_1), \ldots ,\initial_\preceq(f_n)\rangle$. A set of ideal generators for $I$ for which this does hold is called a Gr\"obner basis for $I$ with respect to $\preceq$. A large part of Gr\"obner theory revolves around computing Gr\"obner bases and universal Gr\"obner bases, which are Gr\"obner basis for every possible monomial order on a fixed ideal.
\end{remark}

\begin{ex}
  Let $I = \langle x_2-x_1, x_3 - x_1\rangle \subset \C[x_1, x_2, x_3]$ and let $\preceq$ be the lexicographic order on $\C[x_1, x_2, x_3]$. Then $\initial_\preceq(I) \ne \langle\initial_\preceq(x_2-x_1), \initial_\preceq(x_3-x_1)  \rangle = \langle x_1 \rangle$ since $x_2 - x_1 - (x_3 - x_1) = x_2 - x_3 \in I$ and $\initial_\preceq(x_2-x_3) = x_2 \notin \langle x_1 \rangle$.
\end{ex}

We will now shift from considering a single monomial ordering to considering a plethora of monomial orderings. In order to do this we introduce weight orderings, which allow for a systematic way of considering many, in fact all, possible initial ideals of a certain ideal. The idea is to associate an (in general not total) order on monomials to any vector. 

\noindent Let $R = \C[x_1^{\pm1}, \ldots, x_n^{\pm1}]$ be a Laurent polynomial ring. We make the generalization of working with Laurent polynomial rings instead of regular polynomial rings for technical purposes. However, most ideals we consider in these Laurent polynomial rings will still arise from polynomial rings as being $I\cdot R$ with $I\subset \C[X_1, \ldots, x_n]$. For this reason the reader may find it in their interest to keep regular polynomial rings and ideals in mind.

\begin{defn}[Weight ordering]
    The weight ordering $\preceq_w$ on $R$ associated to $w \in \R^n$ is defined by $x^\alpha \preceq_w x^\beta$ if and only if $\alpha \cdot w \le \beta \cdot w$.
\end{defn}

In this context we call $w$ a weight vector. Now we have to pay slight caution in order to define initial forms, since the order is not total in general and therefore a polynomial does not in general have a unique minimal term. 

\begin{defn}[Initial form w.r.t. a weight vector]
   Let $f= \sum_\alpha c_\alpha x^\alpha \in R$ and $b = \min\limits_{\alpha:c_\alpha \ne0}(w \cdot \alpha)$. The initial form of $f$ with respect to $w$ is defined as:
   \[\initial_w(f) := \sum\limits_{\alpha:w\cdot \alpha = b} c_\alpha x^\alpha\] 
   
\end{defn}

Similarly to monomial orderings, we can define the initial ideal of a polynomial ideal $I \subset R$ as:
\[\initial_w(I) = \langle \initial_w(f) | f \in I \rangle.\] 
We note for future use that the set of all weight vectors $w$ for which there exists a monomial order $\preceq$ such that $\text{in}_\preceq(\text{in}_w(I)) = \initial_\succ(I)$ is called the Gr\"obner region $GR(I)$ of $I$.

\begin{figure}[H]
    \begin{center}      

\begin{tikzpicture}[x=0.7pt,y=0.7pt,yscale=-0.75,xscale=0.75]


\draw  [draw opacity=0][fill={rgb, 255:red, 74; green, 144; blue, 226 }  ,fill opacity=0.55 ] (120.5,479) -- (120.5,299) -- (300.5,299) -- (120,479.5) -- cycle ;
\draw    (300,99.5) -- (300,479.5) ;
\draw    (480,299.5) -- (120,299.5) ;
\draw [color={rgb, 255:red, 74; green, 144; blue, 226 }  ,draw opacity=1 ][line width=2.25]    (300,299.5) -- (120,479.5) ;
\draw [color={rgb, 255:red, 74; green, 144; blue, 226 }  ,draw opacity=1 ][line width=2.25]    (360,479.5) -- (300,299.5) ;
\draw  [draw opacity=0][fill={rgb, 255:red, 74; green, 144; blue, 226 }  ,fill opacity=0.55 ] (120.5,299) -- (350,99.5) -- (350,99.5) -- (350,99.5) -- (301,298.5) -- cycle ;
\draw  [draw opacity=0][fill={rgb, 255:red, 74; green, 144; blue, 226 }  ,fill opacity=0.55 ] (300.5,299) -- (349.5,100) -- (359.5,480) -- (300,299.5) -- cycle ;
\draw  [draw opacity=0][fill={rgb, 255:red, 74; green, 144; blue, 226 }  ,fill opacity=0.55 ] (120.5,479) -- (300,299.5) -- (360,479.5) -- (360,479.5) -- cycle ;
\draw [color={rgb, 255:red, 74; green, 144; blue, 226 }  ,draw opacity=1 ][line width=2.25]    (300,299.5) -- (120,299.5) ;
\draw [color={rgb, 255:red, 74; green, 144; blue, 226 }  ,draw opacity=1 ][line width=2.25]    (350,99.5) -- (300,299.5) ;

\draw (366,272) node [anchor=north west][inner sep=0.75pt]   [align=left] {$\displaystyle \langle y\rangle $};
\draw (246,242) node [anchor=north west][inner sep=0.75pt]   [align=left] {$\displaystyle \langle x^{4} \rangle $};
\draw (176,330) node [anchor=north west][inner sep=0.75pt]   [align=left] {$\displaystyle \langle x^{4} y\rangle $};
\draw (261,378) node [anchor=north west][inner sep=0.75pt]   [align=left] {$\displaystyle \langle x^{3} y^{2} \rangle $};
\draw (345,498) node [anchor=north west][inner sep=0.75pt]   [align=left] {$\displaystyle \langle x^{3} y^{2} +y\rangle $};
\draw (58,488) node [anchor=north west][inner sep=0.75pt]   [align=left] {$\displaystyle \langle x^{4} y+x^{3} y^{2} \rangle $};
\draw (35,288) node [anchor=north west][inner sep=0.75pt]   [align=left] {$\displaystyle \langle x^{4} +x^{4} y\rangle $};
\draw (331,72) node [anchor=north west][inner sep=0.75pt]   [align=left] {$\displaystyle \langle x^{4} +y\rangle $};

\end{tikzpicture}
\caption{Gr\"obner fan of $I = \langle x^4+x^4y-x^3y+x^3y^2+y\rangle$. Cones are labeled by the initial ideals they induce.}
\label{figure:grobnerfan}
    \end{center}
\end{figure}
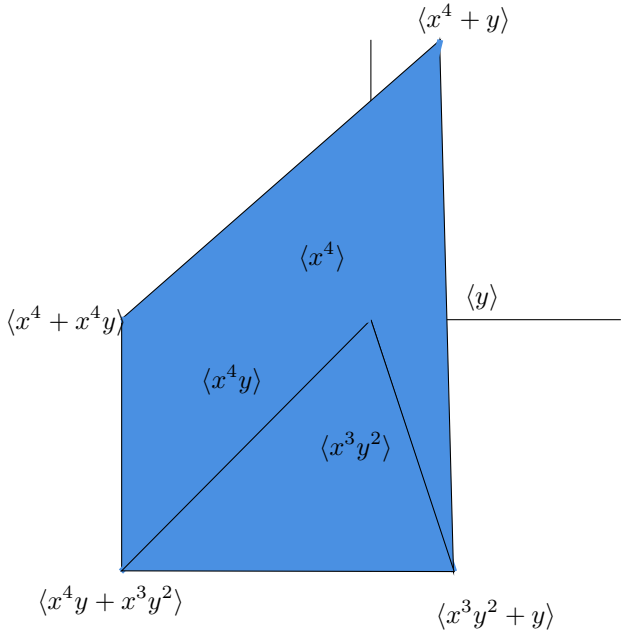

We can now define the Gr\"obner fan, which encodes all initial ideals into a polyhedral fan.
\begin{defn}[Gr\"obner fan]
For any weight vector $w \in \R^n$, the set \[\{v \in \R^n| \initial_v(I) = \initial_w(I)\}\] is an open polyhedral cone. Let $C_w$ be the closure of this cone, then $C_w$ is a closed cone. The set of all $C_w$ for any $w \in \R^n$ forms a polyhedral fan called the Gr\"obner fan of $I$. 
\end{defn}
\begin{ex}
    Consider the principle ideal $I = \langle x^4+x^4y-x^3y+x^3y^2+y\rangle \subset \C[x,y]$. This ideal has four different possible monomial initial ideals:
    \[\initial_{(1,0)}(I) = \langle y \rangle\]
    \[\initial_{(0,1)}(I) = \langle x^4\rangle\]
    \[\initial_{(0,-1)}(I) = \langle x^3y^2\rangle\]
    \[\initial_{(-2,-1)}(I) = \langle x^4y \rangle\]
    Any weight vector in $\R^2\setminus\{0\}$ either yields one of these initial ideals, or a binomial ideal generated by the sum of the generators of two of these monomial ideals. This can be deduced from the Gr\"obner fan of $I$, depicted in \cref{figure:grobnerfan}, which can be computed with a computer algebra system, e.g. Macaulay2. \\
    The Gr\"obner fan contains four two-dimensional cones corresponding to the four monomial ideals, and their intersections form four one dimensional cones, generated by $(1,4), (1,-3), (-1,0)$ and $(-1,-1)$, corresponding to the binomial ideals. Technically the origin is also a cone of this Gr\"obner fan, with $I$ as corresponding initial ideal.\\
    Note how the initial ideal of the intersection of two two-dimensional cones is generated by the sum of the generator of the first two-dimensional cone with the generator of the second two-dimensional cone. This makes sense as the intersection of these two cones is precisely where both of these two monomial generators are minimal with respect to the weight ordering.

\end{ex}
    
\subsubsection{Tropical varieties}

 We can now define tropical varieties, from the viewpoint of Gr\"obner theory.
 In this context, the tropicalization of an ideal $I$ arises as a sub-fan of the Gr\"obner fan of the ideal by considering the set of all weight vectors $w$ such that $\initial_w(I)$ does not contain any monomials.
 
 Let again $R = \C[x_1^{\pm1}, \ldots, x_n^{\pm1}]$ be a Laurent polynomial ring and let $I \subset R$ be an ideal.

\begin{defn}[Tropicalization]
    The tropicalization of $I$ is defined as: 
    \[\trop(I):= \{w\in\R^n\vert\initial_w(I) \text{ does not contain any monomials}\}\]
    
\end{defn}

A priori this does not have a fan structure but it is not difficult to see that it can be endowed with a fan structure analogously to the Gr\"obner fan, and as such is a sub-fan of the Gr\"obner fan. By slight abuse of notation we will denote both the the tropicalization as a set and as a fan by $\trop(I)$.

When working with homogeneous ideals, there are always some weight vectors whose initial ideals are equal to the original ideal. These weight vectors form a subspace called the lineality space.
\begin{defn}[Lineality space]
    The lineality space of an ideal $I$ is the linear subspace $L\subset \R^n$ such that $\initial_w(I) = I$ for all $w\in L$.
\end{defn}

The lineality space is contained in every cone of the tropicalization, and is therefore often quotiented out when studying tropical varieties. 

Tropicalizations are of interest to us because the initial ideals associated to the cones in the tropicalization of an ideal can be of a particularly useful form. In preparation of our discussion of this, we introduce some final terminology.

\begin{defn}[Prime cone]
    A cone $C \subset \trop(I)$ is prime if for any $w \in \relint(C)$, $\initial_w(I)$ is a prime ideal.
\end{defn}

Here $\relint(C)$ denotes the relative interior of $C$, i.e., the interior of $C$ with respect to the euclidean topology of the affine span of $C$. This is a necessary distinction since cones can have lower dimension than the ambient space.

Lastly, we define maximal cones:
\begin{defn}[Maximal cone]
A cone $C_m \subset \trop(I)$ is maximal if its dimension is maximal in $\trop(I)$ i.e. for any cone $C \in \trop(I)$, $\text{dim}(C) \le \text{dim}(C_m)$.
\end{defn}

\begin{ex}
    Consider the Laurent polynomial ring $\C[x^{\pm 1},y^{\pm 1}]$ and the ideal $I = \langle x+xy +y\rangle$. Then $\trop(I)$ contains three maximal cones (of dimension 1), each generated by one of the following vectors: $u =(-1,0), v = (0,-1), w = (1,1)$. Only the cone generated by $w$ is prime, as $\initial_w(I) = \langle x+ y \rangle$ is a prime ideal, while $\initial_u(I) = \langle x+ xy \rangle$ and $\initial_v(I) = \langle xy+ y \rangle$ are clearly not prime ideals.
\tikzset{every picture/.style={line width=0.75pt}}
\begin{figure}
  \centering
  \begin{tikzpicture}[x=0.5pt,y=0.5pt,yscale=-1,xscale=1]

\draw  (250,350) -- (250,50) ;
\draw  (100,200) -- (400,200) ;
\draw [color={rgb, 255:red, 74; green, 93; blue, 226}  ,draw opacity=1 ][line width=2.25]  (250,200) -- (400,50) ;
\draw [color={rgb, 255:red, 255; green, 0; blue, 0 }  ,draw opacity=1 ][line width=2.25]  (250,350) -- (250,200) ;
\draw [color={rgb, 255:red, 255; green, 0; blue, 0 }  ,draw opacity=1 ][line width=2.25]  (100,200) -- (250,200) ;

\end{tikzpicture}
  \caption{Tropicalization of $I= \langle x+xy +y\rangle$ with the prime cone in blue and the non-prime cones in red.}
\end{figure}
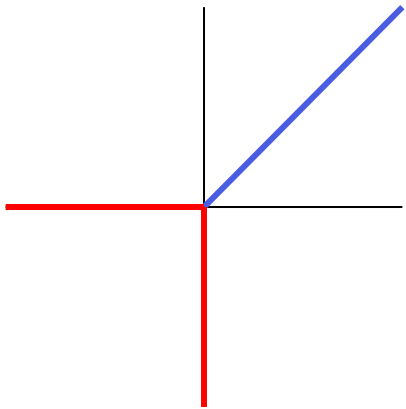
\end{ex}

\begin{remark}
    The general theory of tropical geometry revolves around so-called tropical semi-rings, which replace the usual addition with taking the minimum of the terms and the usual multiplication with regular addition. This leads to the same picture of tropical varieties, but as the Gr\"obner view is more suited for our purposes, we will not elaborate on this. An in-depth treatment of the theory of tropical geometry can be found in \cite{IntroToTropGeom}.
\end{remark}  

\section{Tropical Grassmannians}
In this section we give some background on the tropical Grassmannians, as one of the main focuses of the thesis is studying the tropicalization of $\Gr(3,6)$. Again a good overview of this matter can be found in \cite{IntroToTropGeom}. For more in depth descriptions of higher dimensional tropical Grassmannians, \cite{TropicalGrassmannian}, \cite{TropicalPlanes} are recommended. A more in depth discussion of the space of phylogenetic trees can be found in \cite{PhylogeneticTrees}.

\begin{defn}[Tropical Grassmannians]
    We define the tropical Grassmannian $\trop(\Gr(k,n))$ as the tropicalization of the Pl\"ucker ideal $I_{k,n}$ inside $\C[p_I^{\pm1} : I \in \binom{[n]}{k}]$. 
\end{defn}

The tropical Grassmannians are among the most extensively studied tropical varieties, yielding rich combinatorial descriptions. For example, $\trop(\Gr(2,n))$ has been shown to be equal to a set, known (mostly to biologists) as the space of phylogenetic trees. Nevertheless, higher dimensional tropical Grassmannians are still largely unknown. We give a short description of the space of phylogenetic trees and some insights into the challenges of studying higher dimensional Grassmannians.

\subsection{The tropical Grassmannian trop(Gr(2,n)) and the space of phylogenetic trees}

The Grassmannians with $k=2$ are the simplest Grassmannians (that are not projective spaces), and their tropicalizations are well understood. This is because of a particularly elegant correspondence with the spaces of phylogenetic trees. Let us start with defining a phylogenetic tree.

\begin{defn}[Phylogenetic and trivalent trees]
    A phylogenetic tree on $n$ leaves is a semi labelled tree graph with leaves $1, \ldots, n$ such that each internal node has degree at least three. A phylogenetic tree is trivalent if every internal node has degree exactly three.
\end{defn}

Recall that a graph is a tree if it contains no cycles, and by semi-labelled we mean that only the leaves of the tree are labelled. A phylogenetic tree is commonly used in evolutionary biology, where the leaves would represent certain species, and the internal nodes would represent common ancestors of these species.

\begin{figure}[H]
    \begin{center}

\tikzset{every picture/.style={line width=0.75pt}} 

\begin{tikzpicture}[x=0.75pt,y=0.75pt,yscale=-.7,xscale=.7]

\draw    (300,150) -- (350,200) ;
\draw [shift={(300,150)}, rotate = 45] [color={rgb, 255:red, 0; green, 0; blue, 0 }  ][fill={rgb, 255:red, 0; green, 0; blue, 0 }  ][line width=0.75]      (0, 0) circle [x radius= 3.35, y radius= 3.35]   ;
\draw    (450,250) -- (400,200) ;
\draw [shift={(450,250)}, rotate = 225] [color={rgb, 255:red, 0; green, 0; blue, 0 }  ][fill={rgb, 255:red, 0; green, 0; blue, 0 }  ][line width=0.75]      (0, 0) circle [x radius= 3.35, y radius= 3.35]   ;
\draw    (450,150) -- (400,200) ;
\draw [shift={(450,150)}, rotate = 135] [color={rgb, 255:red, 0; green, 0; blue, 0 }  ][fill={rgb, 255:red, 0; green, 0; blue, 0 }  ][line width=0.75]      (0, 0) circle [x radius= 3.35, y radius= 3.35]   ;
\draw    (300,250) -- (350,200) ;
\draw [shift={(300,250)}, rotate = 315] [color={rgb, 255:red, 0; green, 0; blue, 0 }  ][fill={rgb, 255:red, 0; green, 0; blue, 0 }  ][line width=0.75]      (0, 0) circle [x radius= 3.35, y radius= 3.35]   ;
\draw    (350,200) -- (400,200) ;
\draw    (200,300) -- (250,350) ;
\draw [shift={(200,300)}, rotate = 45] [color={rgb, 255:red, 0; green, 0; blue, 0 }  ][fill={rgb, 255:red, 0; green, 0; blue, 0 }  ][line width=0.75]      (0, 0) circle [x radius= 3.35, y radius= 3.35]   ;
\draw    (400,400) -- (350,350) ;
\draw [shift={(400,400)}, rotate = 225] [color={rgb, 255:red, 0; green, 0; blue, 0 }  ][fill={rgb, 255:red, 0; green, 0; blue, 0 }  ][line width=0.75]      (0, 0) circle [x radius= 3.35, y radius= 3.35]   ;
\draw    (300,400) -- (300,350) ;
\draw [shift={(300,400)}, rotate = 270] [color={rgb, 255:red, 0; green, 0; blue, 0 }  ][fill={rgb, 255:red, 0; green, 0; blue, 0 }  ][line width=0.75]      (0, 0) circle [x radius= 3.35, y radius= 3.35]   ;
\draw    (400,300) -- (350,350) ;
\draw [shift={(400,300)}, rotate = 135] [color={rgb, 255:red, 0; green, 0; blue, 0 }  ][fill={rgb, 255:red, 0; green, 0; blue, 0 }  ][line width=0.75]      (0, 0) circle [x radius= 3.35, y radius= 3.35]   ;
\draw    (200,400) -- (250,350) ;
\draw [shift={(200,400)}, rotate = 315] [color={rgb, 255:red, 0; green, 0; blue, 0 }  ][fill={rgb, 255:red, 0; green, 0; blue, 0 }  ][line width=0.75]      (0, 0) circle [x radius= 3.35, y radius= 3.35]   ;
\draw    (250,350) -- (350,350) ;
\draw    (150,150) -- (200,200) ;
\draw [shift={(150,150)}, rotate = 45] [color={rgb, 255:red, 0; green, 0; blue, 0 }  ][fill={rgb, 255:red, 0; green, 0; blue, 0 }  ][line width=0.75]      (0, 0) circle [x radius= 3.35, y radius= 3.35]   ;
\draw    (150,250) -- (200,200) ;
\draw [shift={(150,250)}, rotate = 315] [color={rgb, 255:red, 0; green, 0; blue, 0 }  ][fill={rgb, 255:red, 0; green, 0; blue, 0 }  ][line width=0.75]      (0, 0) circle [x radius= 3.35, y radius= 3.35]   ;
\draw    (200,200) -- (250,200) ;
\draw [shift={(250,200)}, rotate = 0] [color={rgb, 255:red, 0; green, 0; blue, 0 }  ][fill={rgb, 255:red, 0; green, 0; blue, 0 }  ][line width=0.75]      (0, 0) circle [x radius= 3.35, y radius= 3.35]   ;
        \end{tikzpicture}
    \end{center}
    \caption{Trivalent trees on three, four and five (unlabeled) leaves. Note that there only exists one labeled trivalent tree on three leaves.}
\end{figure}
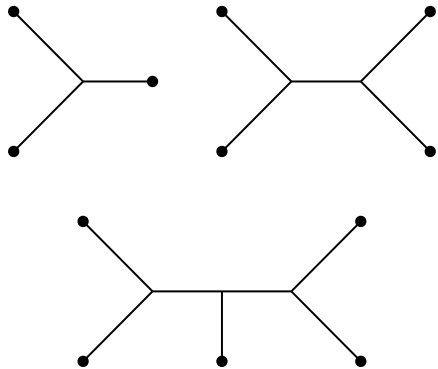

In the biological context, a phylogenetic tree serves mostly as a graphical representation of the evolutionary history of some species, and distance between nodes corresponds to time between evolution. This leads to the key ingredient in the definition of the space of phylogenetic trees on $n$ leaves: the notion of a tree distance.
\begin{defn}[Tree distance]
    A tree distance is a vector in $\R^{n\choose2}$ arising as follows: if $\tau$ is a phylogenetic tree on $n$ leaves, assign to each edge $\epsilon$ of $\tau$ a real number $l_\epsilon$ called its length. A vector $d \in \R^{n\choose2}$ for which $d_{ij}$ is equal to the sum of the lengths of the edges on the path from leaf $i$ to leaf $j$ is a tree distance.
\end{defn}

\begin{remark}
    It is important to note that lengths are not restricted to positive numbers.
\end{remark}

It can be shown that the closure of the set of all tree distances arising from a certain tree forms a polyhedral cone, and the set of all such cones is a polyhedral fan. 

\begin{defn}
    The set of all tree distances arising from trees on $n$ leaves is called the space of phylogenetic trees on $n$ leaves, usually considered with the fan structure as above.
\end{defn}

The maximal cones of this fan correspond to trivalent trees. Intuitively, the reason for this is that any non trivalent tree can be constructed from a trivalent tree by contracting certain edges. Therefore any tree distance obtained from a non trivalent tree can be obtained from the corresponding trivalent trees by setting the lengths of the edges that have to be contracted, equal to zero.

The theorem between $\trop(\Gr(2,n))$ and the space of phylogenetic trees on $n$ leaves can be stated as follows.

\begin{theorem}[\protect{\cite[Theorem 4.3.5]{IntroToTropGeom}}]
    The tropical Grassmannian $\trop(\Gr(2,n))$ is equal to minus the space of phylogenetic trees on $n$ leaves.
\end{theorem}

This correspondence allows study of $\trop(\Gr(2,n))$ by using the concrete trees of the space of phylogenetic trees. Such a concrete description is exceptionally useful when working in the generally abstract field of algebra, and allows for interesting concrete examples of abstract concepts, such as in section 5 of \cite{WallCrossing}.

\subsection{The tropical Grassmannian trop(Gr(3,6)) and higher dimensional Grassmannians}

While efforts have been made to describe higher dimensional Grassmannians in similar ways, complexity increases very quickly. In \cite{TropicalPlanes} the authors introduce arrangements of trees in an effort to find a similar correspondence for $\trop(\Gr(3,n))$. Instead of considering all trees on $n$ leaves, one considers $n$-tuples of trees $(T_1, \ldots, T_n)$ where $T_i$ is a tree on the leaves $[n]\setminus\{i\}$, with compatibility conditions on the metric on these trees. The added complexity is apparent, and while the results of that paper are impressive, they do not yield a complete description of $\trop(\Gr(3,n))$. At the time of writing, research into similar models with $k>3$ is scarce or even nonexistent.

Some of these higher dimensional tropical Grassmannians can still be computed however, through the use of computer algebra systems such as \textit{Macaulay2}\cite{M2}, as they have known symmetries 
. We give some computational results about $\trop(\Gr(3,6))$, which will play an important role in this thesis.

\begin{theorem}[\cite{TropicalGrassmannian}]\label{thm:tropGR(3;6)}
    The tropical Grassmannian $\trop(\Gr(3,6))$ is a 10-dimensional fan in $\R^{20}$ with $1035$ maximal cones, and a 6-dimensional lineality space. 
\end{theorem}

Furthermore these maximal cones can be divided into seven symmetry classes, six of which contain only prime cones and the remaining symmetry class containing only non-prime cones.

\remark{While we have not discussed Grassmannians over general fields, it is important to note that the above results are general in the sense that they do not depend on the characteristic of the field over which these tropical Grassmannians are defined. This only holds for these few simplest cases of Grassmannians however, as it is shown in \cite{TropicalGrassmannian} that $\trop(\Gr(3,7))$ depends on the characteristic of the field over which one is working. This is because of terms with coefficients larger than one appearing in certain initial ideals, resulting in certain initial ideals being monomial in fields with specific finite characteristic, while not containing monomials in fields of other characteristic. This illustrates another challenge in the higher dimensional cases.}

%

%
  \graphicspath{{\chaptersdir/Toric_degenerations/image/}}%
  \chapter{Toric degenerations}\label{ch:Tdegenerations}
In this chapter we introduce toric ideals and toric degenerations. We then show how they can be obtained from tropical varieties and from Newton-Okounkov bodies. Finally, we state the link between these two methods and state the wall-crossing theorem for prime cones.

\section{Basic toric geometry}

We start with some background on toric ideals. This is a short summary of relevant definitions, propositions and theorems of chapter four of \cite{GrobnerBasisAndConvexPolytopes}, which we refer to for proofs of the statements we give.

\begin{defn}[Toric ideal]
Given a $d\times n$ matrix $A \in \Z^{d\times n}$, we construct a map from the polynomial ring $\C[x_1, ... ,x_n]$ to the Laurent polynomial ring $\C[t_1^\pm ,..., t_d^\pm]$ :
\[\pi: \C[x_1, ... ,x_n] \to \C[t_1^\pm ,..., t_d^\pm]\quad
x^u \mapsto t^{Au}\]
    The toric ideal $I_A$ of $A$ is the kernel of the map $\pi$.
\end{defn}

\remark{We will only really concern ourselves with homogeneous toric ideals, as we want to work over projective space. For this, it suffices to require that $A$ contains a row that is identically one.}

A toric ideal $I_A$ induces an affine toric variety $V(I_A)$. These ideals and varieties are studied because of their combinatorial descriptions, some of which we specify below. Note that one has to be careful when reading about toric varieties, as some authors define their toric varieties to have an additional property called normality. This is usually done when studying the varieties in a more algebraic geometry oriented context. Since we focus more on the combinatorial side, we stick with the above definition and therefore do not assume normality. For a discussion on normality and the relation between the two definitions we refer to chapter 14 of \cite{GrobnerBasisAndConvexPolytopes}.

We state some lemmas characterizing the form of toric ideals as being generated by binomials.

\begin{lemma}[\cite{GrobnerBasisAndConvexPolytopes} Lemma 4.1-4.4]
    The toric ideal $I_A$ is spanned as a $\C$-vector space by the set of binomials
    \[\{x^u-x^v ; u,v\in \N^n \text{ s.t. } Au=Av\}.\] Furthermore, for any $u\in \Z^n$, we can write $u = u^-u^-$ where $u^+=\frac{1}{2}(u+|u|)$ and $u^-=-\frac{1}{2}(u-|u|)$, and there exists a finite set $G_A \subset \ker(A)$ (where $A$ is seen as a map from $\Z^n$ to $\Z^d$) such that\[I_A = \langle x^{u^+}-x^{u^-}\vert u \in G_A \rangle.\]
\end{lemma}


One of the most important theorems concerning toric ideals gives an explicit form of the Hilbert polynomial of $V(I_A)$ in terms of the polytope $Q = \mathrm{conv}(A)$. In the following we assume $A = Q \cap \Z^d$. 

\begin{defn}[Ehrhart polynomial] 
    Assuming $I_A$ is homogeneous, we define a function \[E_Q: \N \to \N \quad
    r \to | \Z A \cap r\cdot Q|\] where $\Z A$ is the integral span of the columns of $A$ (the lattice generated by the columns of $A$) and $r\cdot Q$ is the $r$-dilation of $Q$, given by the set$\{r\cdot v \vert v\in Q\}$. It can be shown that this is a polynomial, which we call the Ehrhart polynomial of $Q$, and we let $q$ denote its degree.
\end{defn}

\remark{The Ehrhart polynomial can also be defined for non-lattice polytopes. However, it will no longer be polynomial. In this context it is called the Ehrhart function.}

We need a further condition on the polytope $Q$.
\begin{defn}[Normal polytope]
    A lattice polytope $P\subset \R^n$ is normal if for every $k\in \N$, all lattice points of $k\cdot P$ can be written as the sum of $k$ lattice points of $P$.
\end{defn}

We can now state the theorem yielding an explicit form for the Hilbert polynomial of a toric ideal.

\begin{theorem}[\cite{GrobnerBasisAndConvexPolytopes} Theorem 4.16]
    Assume $Q$ is normal. The Hilbert polynomial of the variety $V(I_A)$ is equal to the Ehrhart polynomial of $Q$, and in particular, the degree of $V(I_A)$ is equal to the so-called normalized volume of $Q$, which is $q!$ times the leading coefficient of $E_Q$.
\end{theorem}

\section{Degenerations}

In this section we define degenerations and show how Gr\"obner theory can be used to obtain degenerations of an ideal to any of its initial ideals. In particular, if the initial ideal is toric, a toric degeneration is obtained.

We begin by defining a degeneration.

\begin{defn}[Degeneration]
    A degeneration of an ideal $I$ to an ideal $I_0$ is a flat family of ideals $\tilde{I}_t$ over $\mathbb{A}^1$ such that the generic fibers, $t\ne0$, are isomorphic to $I$ and the special fiber, $t=0$, is equal to $I_0$.
\end{defn}

Here a flat family is in the scheme theoretic sense, i.e. we assume there is a flat morphism of schemes $\mathrm{Proj}(\C[x_0,\ldots,x_n])\to \mathbb{A}^1$ such that the generic fibers are isomorphic to $\mathrm{Proj}(\C[x_0,\ldots,x_n]/I)$ and the special fiber is $\mathrm{Proj}(\C[x_0,\ldots,x_n]/I_0)$, where $I,I_0$ are prime and homogeneous. For a more in depth treatment of flat families we refer to \cite[chapter III section 9]{Hartshorne}.

The idea of a degeneration is to smoothly deform a variety into a variety with certain desirable properties, in such a way that invariants of the original variety are preserved. The flatness of the family is what guarantees preservation of certain invariants, such as degree, Hilbert polynomial, etc. Note that this means that existence of a degeneration is not always guaranteed, as one can not, for example, degenerate an ideal to an ideal of different degree.

\subsection{Gr\"obner degenerations}

One way to obtain degenerations of an ideal $I$ is by studying its Gr\"obner fan. Recall that each cone $C$ of the Gr\"obner fan of $I$ has an associated initial ideal $\initial_C(I)$ of $I$, which is obtained by taking a weight vector $w\in \relint{C}$ and calculating $\initial_C(I) = \initial_w(I)$. Gr\"obner theory yields an explicit construction for a degeneration of $I$ to $\initial_C(I)$ for any cone in the Gr\"obner fan, as shown in the following theorem.
\begin{theorem}[\cite{eisenbud1995commutative} Theorem 15.17]
    Let $I$ be an ideal, $C$ a cone in the Gr\"obner fan of $I$ and $w \in \relint(C)$. The family of ideals \[\tilde{I}_t:= \left\langle t^{-\min_u(w\cdot u)}f(t^{w_1}x_1, ..., t^{w_n}x_n)\vert f = \sum a_ux^u \in I \right\rangle\] is a degeneration of $I$ to $\initial_C(I)$.
\end{theorem}

We are particularly interested in finding degenerations to toric ideals, i.e. toric degenerations. For this purpose, we can immediately disregard cones whose initial ideals are monomial, as these can not yield toric degenerations. 

This leaves us with the subfan of the Gr\"obner fan given by all cones whose initial ideals do not contain monomials, which we know to be the tropicalization of our ideal. This motivates our interest in tropicalizations of ideals: the cones in the tropicalization of an ideal yield candidates for possible toric degenerations of that ideal. We are particularly interested in the prime, maximal cones of the tropical varieties, as these will precisely be the cones yielding toric degenerations through Gr\"obner theory. Finding these prime, maximal cones is not always an easy task however.\
\begin{ex}
    Consider $I = \langle xy+yz-zw\rangle \subset \C[x,y,z,w]$. Note that the weight vector $u = (0,1,1,0)$ yields the toric initial ideal $\initial_u(I) = \langle xy-zw\rangle$. Using the above theorem, we can construct a toric degeneration of $I$ to $\initial_u(I) = \langle xy-zw\rangle$ as:
    \[\tilde{I}_t= \left\langle t^{-1}(t^0xt^1y+ t^1yt^1z- t^1zt^0w)\right\rangle =\left\langle xy+ tyz- zw\right\rangle\]
\end{ex}

\section{Valuations and Newton-Okounkov bodies}
In order to find prime maximal cones, we will exploit a connection between the theory of tropical varieties and the theory of Newton-Okounkov bodies. We start by introducing the relevant concepts surrounding Newton-Okounkov bodies, and then show the link between cones in the tropical variety of an ideal $I$ and Newton-Okounkov polytopes of the coordinate ring of $V(I)$. We closely follow \cite{Kaveh2009NewtonOkounkovBS} and \cite{KhovanskiiBases} for this section.

In the theory of Newton-Okounkov bodies one studies valuations on algebras, which are maps of the algebra to some discrete set. The images of these valuations will form semigroups, and will give rise to convex objects, Newton-Okounkov bodies. The valuations of interest will be those giving rise to finitely generated semigroups, as their Newton-Okounkov bodies will be convex polytopes. When this is the case for a valuation on a coordinate ring $\C[x_0,...,x_n]/I$, it will yield a toric degeneration of $V(I)$ to the toric variety associated to the Newton-Okounkov polytope.

\subsection{Valuations and filtrations}
Let $A$ be a (finitely generated) algebra of Krull dimension $d$, and let $\prec$ be a fixed total order on $\Q^r$ with $r \in \N$.
\begin{defn}[Valuation]
    A valuation over $\C$ is a function $\nu: A \setminus \{0\} \to \mathbb{Q}^r$ such that for all $f,g\in A \setminus\{0\}$, the following properties hold:
    \begin{enumerate}
	\item  if $f+g \ne 0$, then $\nu(f+g) \succeq \min\{\nu(f),\nu(g)\}$
        \item $\nu(f \cdot g) = \nu(f) + \nu(g)$
        \item $\forall c \in \C\setminus \{0\}$: $\nu(cf) = \nu(f)$ or, equivalently, $\nu(c) = 0$
    \end{enumerate}
\end{defn}

We will assume all valuations to be discrete, i.e. their image is a discrete subset of $\Q^r$.
A valuation has an associated multiplicative filtration $\mathcal{F}_\nu$ on $A$ defined, for every $a \in \Q^r$, as follows: 
\[F_{\nu\succeq a} = \{f \in A | \nu(f)\succeq a\}\cup\{0\}\text{  and  } F_{\nu\succ a} = \{f \in A | \nu(f)\succ a\}\cup\{0\} \]
A valuation $\nu$ is said to have one-dimensional leaves if $F_{\nu \succeq a}/F_{\nu\succ a}$ is at most one dimensional as a vector space. Through its filtration, a valuation then leads to a graded algebra structure on $A$: $gr_\nu(A) = \bigoplus\limits_{a \in \Q^r}F_{\nu \succeq a}/F_{\nu\succ a}$.
\begin{remark}
	One way to think about this construction is as a generalized version of Gr\"obner theory. A polynomial ideal $I$ becomes  an algebra $A$, a monomial ordering $\preceq$ on $I$ becomes a valuation $\nu$ and the initial ideal $\initial(I)$ becomes the graded algebra $gr_\nu(A)$.
\end{remark}

Let us look at some examples:
\begin{ex}
    Consider $ \C[x_1, ..., x_r]$ as an algebra, and let $\succeq$ be a fixed order on $\Q^r$. We can define a valuation as follows: let $f =\sum_{\alpha} c_\alpha x^\alpha\in \C[x_1, ..., x_n]$, define $\nu: A \to \Q^r: f \mapsto \min_{\succeq}(\alpha:c_\alpha \ne 0)$. This valuation has one-dimensional leaves.
\end{ex}

\begin{ex}\label{HomVal}
    Again consider $ \C[x_1, ..., x_r]$ as an algebra. Let $\succeq_m$ be the monomial ordering induced by our fixed order on $\Q^r$, i.e. $x^\alpha  \succeq_m x^\beta$ if and only if $\alpha \succeq  \beta$.  We can define a valuation $\nu$ as follows. Let $f \in \C[x_1, ..., x_n]$ and let $\initial_{\succeq_m}(f) = c_\alpha x^\alpha$, define $\nu(f) = -\alpha$. 
    This valuation again has one dimensional leaves. Notice as well that this valuation essentially is a reverse order comparison  of the degrees of polynomials: if $\left(\text{in}_{\succeq_m}(f_1)\right) \prec_m \left(\text{in}_{\succeq_m}(f_2)\right)$ then $\nu(f_1) \succ \nu(f_2)$.
\end{ex}

The valuations that we consider will always be of a similar form to the valuation in \cref{HomVal}. Concretely, we discuss valuations on positively graded algebras which are homogeneous with respect to the degree:
\begin{defn}[Homogeneous valuation]
    Let $A = \bigoplus\limits_{n\in \N}A_n$ be a positively graded algebra, then a valuation $\nu$ on $A$ is homogeneous with respect to the degree if the following holds: for any $f_1,f_2 \in A\setminus\{0\}$, if $\text{deg}(f_1) < \text{deg}(f_2)$, then $\nu(f_1) \succ \nu(f_2)$.
\end{defn}

Specifically, we will consider valuations to $\Q^r$ of the form \begin{equation}\label{Valuationform}
	\nu(f) = (\text{deg}(f),\nu'(f))
\end{equation}
where $\nu'$ is a general (not necessarily homogeneous) valuation to $\Q^{r-1}$.  The total order on $\Q^r$ will be taken to be a modified version of the lexicographic order, that is reversed for the first coordinate and the same for all others. This yields a valuation that first compares the degree of two elements, with a reversal in the order as above (i.e. if $\text{deg}(f_1) < \text{deg}(f_2)$ then $\nu(f_1) \succ \nu(f_2)$), and then breaks ties with the other coordinates. 

\subsection{Newton-Okounkov bodies}

We have introduced valuations on algebras, and now continue to the convex objects they give rise to: Newton-Okounkov cones and bodies. 

We start of by considering the images of such valuations. The image of a (discrete) valuation $\nu$ on an algebra $A$ is a (discrete) additive semigroup $S(A,\nu) := \nu(A\setminus\{0\})$ which sits inside $\Q^r$.

\begin{defn}[Semigroup]
    A semigroup is a set $S$ equipped with an associative binary operation $\cdot: S\times S \to S$.
\end{defn}
A semigroup is essentially a group without guaranteed existence of inverses or an identity element. Semigroups are usually studied because they give rise to lattices, which arise in many areas of mathematics. Notably, toric varieties also have connections to semigroups.
\begin{remark}
    There are names for each of the structures obtained from a group without requiring any subset of the three properties we ask of a group operation: associativity, inverses and an identity element. For example, the binary operation on a magma requires none of these properties, while a monoid is a semigroup with an identity element.
\end{remark}
A semigroup can be turned into a group by adding all inverses and an identity element. As such, it inherits a notion of rank as the rank of its so-called group completion. This allows us to define the rank of a valuation.

\begin{defn}[Rank of a valuation]
    The rank of a valuation $\nu$ on an algebra $A$ is the rank of the group generated by its value semigroup $S(A,\nu)$. A valuation of rank $d = \text{dim}(A)$ is said to have full-rank.
\end{defn}

Full-rank valuations are particularly desirable, especially when working over an algebraically closed field and assuming $A$ is a domain. This is because, in that case, two results from \cite{KhovanskiiBases} show that a full rank valuation has one-dimensional leaves, and having one-dimensional leaves implies $gr_\nu(A)$
is isomorphic to the semigroup algebra generated by $S$. A semigroup algebra is a construction equivalent to the group algebra for a group. Notably, some semigroups can give rise to a toric variety as the spectrum of the associated semigroup algebra. This is, intuitively, how valuations on a coordinate ring $A \cong \C[x_1,.., x_n]/I$ yield toric degenerations of $I$. To see which valuations yield toric degenerations, we study associated convex bodies called the Newton-Okounkov cone and the Newton-Okounkov body of a valuation on $A$.

\begin{defn}[Newton-Okounkov cone]
    The Newton-Okounkov cone associated to $A$ and $\nu$ is the cone generated by the set $S(A,\nu) \subset \R^r$, i.e. 
\[\text{Cone}(S(A,\nu)) = \left\{\sum_i t_is_i|s_i\in S(A,\nu), t_i \in \R_{\ge 0}\right\}.\]
\end{defn}

\begin{defn}[Newton-Okounkov body]
    The Newton-Okounkov body associated to $S(A,\nu)$ is defined as: \[\Delta(A,\nu) :=  \text{Cone}(S(A,\nu)) \cap \{x_1 = 1\}\]
\end{defn}
\begin{remark}
We will assume $A$ is contained in a finitely generated algebra, which guarantees that $\Delta(A,\nu)$ is bounded according to \cite{KhovanskiiBases}.
 \end{remark}
 
 It is of fundamental importance to note that the Newton-Okounkov body is a convex \textit{body}: it is not always a polytope. The same goes for the Newton-Okounkov cone, it is always convex but in general not polyhedral. The problem lies in the aspect of finiteness we require of polyhedra, i.e. that they are an intersection of a finite amount of half-spaces. As such it makes sense that when $S(A,\nu)$ is finitely generated, then $\text{Cone}(S(A,\nu))$ is polyhedral and $\Delta(A,\nu)$ is a polytope. However, notably, and perhaps counter intuitively, the converse does not hold, as shown in  \cite{ObodiesOfFGDivisors} and \cite{Seppanen2014OkounkovBF} for example. 
\begin{ex}\label{ex : valuation}
  Consider $A \cong \C[x,y,z]/I$ where $I = \langle  x^2+xy+xz+z^2\rangle$. Consider the matrix \[M :=\begin{bmatrix}
    1 & 1& 1\\
    2 & 0  &  1
  \end{bmatrix}\]
   let $f =\sum_{\alpha} c_\alpha x^\alpha\in \C[x,y,z]$, we define a valuation on $\C[x,y,z]$ as follows: 
   \[\tilde{\nu}_M(f) := \min\{M\alpha|c_\alpha \ne 0\}\]
   where the minimum is taken with respect to lexicographic order on $\Q^2$ that is reversed for the first coordinate. In the next section we will see that this valuation on $\C[x,y,z]$ induces a valuation $\nu_M$ on $A$, whose value semigroup is generated by the column vectors of $M$. The Newton-Okounkov cone associated to $A$ and $\nu_M$ is then the   cone generated by $(1,0)$ and $(1,2)$, and the Newton-Okounkov body $\Delta(A,\nu_M)$ is the line segment from  $(1,0)$ to $(1,2)$ as depicted in figure 2.
   \begin{figure}[h]
       \centering
       \begin{tikzpicture}[x=0.5pt,y=0.5pt,yscale=-0.75,xscale=0.75]
\draw    (250,300) -- (250,0) ;
\draw    (50,200) -- (450,200) ;
\draw  [draw opacity=0][fill={rgb, 255:red, 74; green, 144; blue, 226 }  ,fill opacity=0.41 ] (450,100) -- (250,200) -- (250,0) -- (450,0) -- cycle ;
\draw [color={rgb, 255:red, 74; green, 144; blue, 226 }  ,draw opacity=1 ][line width=2.25]    (250,0) -- (250,200) ;
\draw [color={rgb, 255:red, 74; green, 144; blue, 226 }  ,draw opacity=1 ][line width=2.25]    (450,100) -- (250,200) ;
\draw [color={rgb, 255:red, 252; green, 0; blue, 0 }  ,draw opacity=1 ][line width=1.5]    (250,150) -- (350,150) ;
\draw [shift={(350,150)}, rotate = 0] [color={rgb, 255:red, 252; green, 0; blue, 0 }  ,draw opacity=1 ][fill={rgb, 255:red, 252; green, 0; blue, 0 }  ,fill opacity=1 ][line width=1.5]      (0, 0) circle [x radius= 4.36, y radius= 4.36]   ;
\draw [shift={(250,150)}, rotate = 0] [color={rgb, 255:red, 252; green, 0; blue, 0 }  ,draw opacity=1 ][fill={rgb, 255:red, 252; green, 0; blue, 0 }  ,fill opacity=1 ][line width=1.5]      (0, 0) circle [x radius= 4.36, y radius= 4.36]   ;

\end{tikzpicture}
\caption{The Newton-Okounkov cone (blue) and body (red) of example 6.}
   \end{figure}
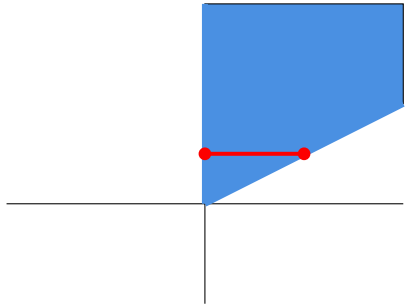
\end{ex}

We now define Khovanskii bases, also referred to as SAGBI bases, which stands for Subalgebra Analogue of Gr\"obner
Basis for Ideals. 

\begin{defn}[Khovanskii basis]
	A subset $\mathcal{B} \subset A\setminus\{0\}$ is a Khovanskii basis for $(A,\nu)$ if the image of $\mathcal{B}$ in $gr_\nu(A)$ generates it as an algebra.
\end{defn}

Note that this is a stronger condition than that of being a Gr\"obner basis, since a generating set for an ideal does not necessarily generate that ideal as an algebra. Khovanskii bases are not assumed to be finite, but the case of existence of a finite Khovanskii basis is what interests us, as we have the following result from \cite{KhovanskiiBases}:
\begin{lemma}[\cite{KhovanskiiBases} Lemma 2.10] Let $\mathcal{B}$ be a Khovanskii basis for $(A, \nu)$. Then the set of values $\{\nu(b)|b\in\mathcal{B}\}$ generates $S(A,\nu)$ as a semigroup.
\end{lemma}

It follows that existence of a finite Khovanskii basis implies that $\Delta(A,\nu)$ is a polytope. We refer to \cite{KhovanskiiBases} section 2.2. for necessary and sufficient conditions for a set to be a Khovanskii basis, and for an algorithm to find a Khovanskii basis. The final notion of this section is that of an adapted basis:

\begin{defn}[Adapted basis]
	A vector space basis $\mathcal{B}$ for $A$ is an adapted basis for $(A,\nu)$ if its image in $gr_\nu(A)$ is a vector space basis of $gr_\nu(A)$.
\end{defn}

\begin{remark}
    With the Khovanskii bases taking the role of Gr\"obner bases, adapted bases can be thought of as the analogue for the standard monomial bases arising in Gr\"obner theory.
\end{remark}

\subsection{Quasi-valuations and weight valuations}
\label{Quasi-valuation}

To link Newton-Okounkov bodies with the tropical picture we need the following slight generalization of valuations:

\begin{defn}[Quasivaluation]
	A quasivaluation is a function $\nu: A \setminus \{0\} \to \mathbb{Q}^r\cup\{\infty\}$ such that $\forall f,g\in A \setminus\{0\}$:
	\begin{enumerate}
		\item  if $f+g \ne 0$, then $\nu(f+g) \succeq \min\{\nu(f),\nu(g)\}$
		\item $\nu(f \cdot g) \succeq \nu(f) + \nu(g)$
		\item $\forall c \in \C\setminus \{0\}$: $\nu(cf) = \nu(f)$ or equivalently $\nu(c) = 0$
	\end{enumerate}
\end{defn}

Note that the only difference with a valuation is in the second property where we now allow the valuation of a product to be greater than the sum of the valuations of the terms.\\
A quasi-valuations again gives rise to a filtration on $A$ in the same way as a regular valuation. It is also possible to go from a filtration to a quasi-valuation:
for a decreasing algebra filtration $\mathcal{F} := \{F_a\}_{a \in \Q^r} $ we can define \[\nu_{\mathcal{F}}(f) := \max\{a \in \Q^r|f \in F_a\}\] which can be shown to be a quasi-valuation.\\

We now detail a construction of quasi-valuations that is central to the link with tropical geometry.

Let $\C[x_1,...,x_n]/I\cong A$ be a fixed presentation for $A$. Fix a total ordering $\succ$ on $\Q^r$. 
For a matrix $M \in \Q^{r\times n}$ we define a valuation $\tilde{\nu}_M$ on $C[x_1, ..., x_n]$ as follows: let \[f =\sum_{\alpha\in A} c_\alpha x^\alpha\in \C[x_1, ..., x_n]\] then we define the weight valuation induced by $M$ as: \[\tilde{\nu}_M(f) := \min\{M\alpha|c_\alpha \ne 0\}\]
where the minimum is taken with respect to $\succ$. Similarly to \cref{Gr\"obnertheory}, define the initial form of $f$ with respect to $M$ to be \[\initial_M(f) = \sum\limits_{\alpha:M\alpha = \tilde{\nu}_M(f)} c_\alpha x^\alpha\] and the initial ideal of $I$ with respect to $M$ to be the ideal \[\initial_M(I) = \langle \initial_M(f)|f\in I\rangle.\]
Denote $\tilde{\mathcal{F}}_M$ for the filtration on $\C[x_1, ..., x_n]$ induced by $\tilde{\nu}_M$. We can project the filtration $\tilde{\mathcal{F}}_M$ to get a filtration on $A$ via our chosen representation, we call this the weight filtration associated to $M$ on $A$. Finally, this weight filtration gives us the weight quasi-valuation $\nu_M$ on $A$, which, notably, is not always a valuation.
\section{Wall-crossing for prime cones}

Now that we have introduced the relevant concepts, we can introduce the theorems on which our research hinges. The theorems provide a link between tropical geometry and the theory of Newton-Okounkov bodies, and show how that link can be used to wall-cross between a prime cones.

\subsection{Valuations from prime cones}

The first result is from \cite{KhovanskiiBases} and it provides the fundamental link between tropical geometry and Newton-Okounkov bodies. 

\begin{theorem}[\cite{KhovanskiiBases} Propositions 4.2 and 4.8]\label{TropNOlink} Let $C \subset \trop(I)$ be a maximal prime cone of $\trop(I)$. Let $\{u_1, ...,u_{d}\}\subset C$ be a collection of rational vectors which span the span of $C$. Let $M \in Q^{d \times n} $ be the matrix whose row vectors are $u_i$. Let $\nu_M: A \setminus\{0\}\to \Q^{d} $ denote the corresponding weight quasi-valuation. Then $\nu_M$ is a valuation and the following hold:
	\begin{enumerate}
		\item $gr_{\nu_M}(A) \cong \C[x_1, ..., x_n]/\text{in}_M(I)$ as $\Q^{d}$-graded algebras
		\item The column vectors of $M$ generate the value semigroup $S(A, \nu_M)$
		\item If $C$ lies in $GR(I)$, the projection via $\pi$ of $\mathcal{S}(<,I)$ is an adapted basis for $\nu_M$ where $C_<$ is a max-dimensional cone in the Gr\"obner fan of $I$ with $C \subset C_<$.
	\end{enumerate}
\end{theorem}
This result allows us to construct valuations from prime cones, and hence construct Newton-Okounkov bodies from prime cones. Section 5 of \cite{KhovanskiiBases} provides a similar result for the reverse direction, constructing a prime cone from a valuation, but this is of less relevance to us. It is important to note that when working with homogeneous ideals (as will be the case for us), the Gr\"obner region of $I$ is all of $\Q^n$ (see \cite{GrobnerBasisAndConvexPolytopes} proposition 1.12) and the condition in $3.$ is guaranteed to be satisfied.

\begin{ex}
    Consider again the case of \cref{ex : valuation}: $A \cong \C[x,y,z]/I$ where $I = \langle x^2+xy+xz+z^2\rangle$. In \cref{ex : valuation} we stated that the matrix \[M :=\begin{bmatrix}
    1 & 1& 1\\
    2 & 0  &  1
  \end{bmatrix}\] and its accompanying weight valuation $\tilde{\nu}_M$ on $\C[x,y,z]$ induce a valuation $\nu_M$ on $A$. Here $\nu_M$ is in fact the weight quasi-valuation defined in \cref{Quasi-valuation}. We show that this is indeed a valuation by using \cref{TropNOlink}: the tropicalization $\trop(I)$ consists of the three cones \[C_1 = \R(1,1,1)+\R_{\ge0}(2,0,1) = \R(1,1,1)+\R_{\ge0}(0,-2,-1)\] \[C_2 = \R(1,1,1)+\R_{\ge0}(0,1,0)\] \[C_3 =\R(1,1,1)+\R_{\ge0}(0,0,1)\] which are depicted in figure 3. Note that $C_1$ is prime as it gives the prime initial ideal $\langle xy+z^2\rangle$. Now note that the row vectors of $M$ are the vectors $(1,1,1)$ and $(2,0,1)$, which span the span of $C_2$. Theorem 1 now exactly gives us that $\nu_M$ is a valuation and that $S(A,\nu_M)$ is generated by the column vectors of $M$.
  \begin{figure}[h]
      \centering
      \begin{tikzpicture}[x=0.5pt,y=0.5pt,yscale=-0.70,xscale=0.70]

\draw    (250,400) -- (250,0) ;
\draw    (50,200) -- (450,200) ;
\draw [color={rgb, 255:red, 74; green, 144; blue, 226 }  ,draw opacity=1 ][line width=2.25]    (450,200) -- (250,200) ;
\draw [color={rgb, 255:red, 74; green, 144; blue, 226 }  ,draw opacity=1 ][line width=2.25]    (250,200) -- (250,0) ;
\draw [color={rgb, 255:red, 74; green, 144; blue, 226 }  ,draw opacity=1 ][line width=2.25]    (250,200) -- (50,300) ;
\end{tikzpicture}
\caption{The tropicalization (with lineality space quotiented out) $\trop(I)/\R(1,1,1)$ depicted in the plane spanned by $(0,1,0)$ and $(0,0,1)$.}
  \end{figure}
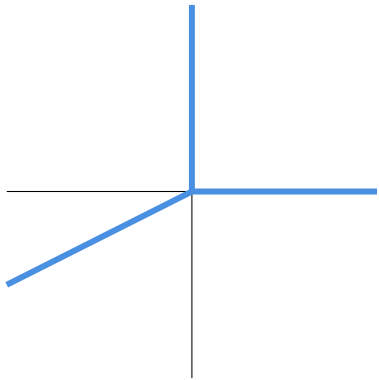
\end{ex}

\subsection{The Wall-crossing theorem}

The above results were the main motivation for the wall-crossing paper \cite{WallCrossing} by Escobar and Harada. The relation described in \cite{WallCrossing} starts from the construction that assigns a valuation, and therefore a Newton-Okounkov body, on $A = \C[x_1, ..., x_n]/I$ to a cone in $\trop(I)$. It is then shown that for two adjacent maximal prime cones in $\trop(I)$, there are two piece-wise linear maps between the Newton-Okounkov polytopes associated to the adjacent prime maximal cones.\\
\\
Let $C_1,C_2 \subset \trop(I)$ be two distinct maximal cones in the tropicalization of $I$.
\begin{defn}
	The cones $C_1,C_2$ are called \textit{adjacent} if they share a face of co-dimension 1.
\end{defn}
Let us now take $C_1$ and $C_2$ to be adjacent prime maximal cones, and denote their common face by $C = C_1 \cap C_2$. By taking $d-1$ vectors in $C$ such that they span the linear span of $C$, we can perform the construction of the weight matrix in \cref{TropNOlink} for both $C_1$ and $C_2$ such that their matrices differ only in the last row.

More explicitly: take $u_1, ..., u_{d-1}$ to be the chosen vectors in $C$. Here we choose $u_1 = (1, ..., 1)$, which is possible because $I$ is homogeneous. Then take $w_1$ (respectively $w_2$) such that together with $\{u_1, ..., u_{d-1}\}$ they span the linear span of $C_1$ (respectively $C_2$).

Now let $M$ be the $(d-1)$ by $n$ matrix whose row vectors are $u_1, ..., u_{d-1}$ and let $M_1$ (respectively $M_2$) be the $d$ by $n$ matrix obtained by adding $w_1$ (respectively $w_2$) as a new row at the bottom of $M$.\\
Denote $\nu_M, \nu_{M_1}, \nu_{M_2}$ for the induced weight quasi-valuations on $A$, then \cref{TropNOlink} implies that $\nu_{M_1}$ and $\nu_{M_2}$ are valuations. Moreover, by taking $u_1 = (1, ..., 1)$, we have ensured that $\nu_{M_1}$ and $\nu_{M_2}$ are of the form described in \cref{Valuationform}, and that the first coordinate of the column vectors of $M_1$ and $M_2$ is always equal to one. 

From the second point of \cref{TropNOlink} it then follows that the Newton-Okounkov body $\Delta(A,\nu_{M_1})$ (respectively $\Delta(A,\nu_{M_2})$) is equal to the convex hull of the column vectors of $M_1$ (respectively $M_2$). It follows from the construction that the projections onto the first $d-1$ coordinates, which we denote by $p_{[1,d-1]}$, of these Newton-Okounkov bodies are both equal to the convex hull of the column vectors of $M$. 
\begin{remark}
    It is actually possible to show that $\nu_M$ is also a valuation, albeit not of full rank, and its Newton-Okounkov body $\Delta(A,\nu_M)$ is also equal to the convex hull of the column vectors of $M$.
\end{remark} 

We now state the wall-crossing theorem, which is the main theorem of \cite{WallCrossing} and the main motivation for our question of wall-crossing for non-prime cones:
\begin{theorem}[\cite{WallCrossing} Theorem 2.7]\label{thm : wall-crossing}
Let $A = \bigoplus_kA_k$ be a positively graded algebra over $\C$, and assume $A$ is an integral domain and has Krull dimension $d$. Let $\mathcal{B} = \{b_1, . . . , b_n\}$ be a subset of $A_1$ (the homogeneous degree 1 elements of $A$) which generate $A$ as an algebra. Let $I$ be the homogeneous ideal such that the presentation induced by $\mathcal{B}$ is $A \cong \C[x_1, . . . , x_n]/I$, and let $\trop(I)$ denote its tropicalization. Suppose that $C_1$ and $C_2$ are two maximal	prime cones in $\trop(I)$ that share a face $C$ of codimension 1. Let $M_1$, $M_2$, and $M$ be the matrices described above and $\nu_{M_1}$ , $\nu_{M_2}$ and $\nu_{M}$ the corresponding weight valuations on $A$. Let $\Delta(A, \nu_{M_1} ) \subset \{1\} \times \R^{d-1}, \Delta(A, \nu_{M_2} ) \subset \{1\} \times \R^{d-1}$ and $\Delta(A, \nu_{M} ) \subset \{1\} \times \R^{d-2}$ denote the corresponding Newton-Okounkov bodies. Let $p_{[1,d-1]} : \R^{d} \to \R^{d-1}$
	denote the linear projection obtained by deleting the last coordinate. Then \[p_{[1,d-1]}(\Delta(A, \nu_{M_1} )) = p_{[1,d-1]}(\Delta(A, \nu_{M_2} )) = \Delta(A, \nu_{M} )\]
	and for any $\xi \in \Delta(A, \nu_{M})$, the Euclidean lengths of the fibers $p_{[1,d-1]}^{-1}(\xi)\cap \Delta(A, \nu_{M_1})$ and $p_{[1,d-1]}^{-1}(\xi)\cap \Delta(A, \nu_{M_2})$ are
	equal, up to a global constant which is independent of $\xi$. Moreover, there exist two piece-wise linear identifications
	$S_{12}: \Delta(A, \nu_{M_1}) \to \Delta(A, \nu_{M_2}) $ and $F_{12}: \Delta(A, \nu_{M_1}) \to \Delta(A, \nu_{M_2})$ , called the “shift map” and the “flip map” respectively, which have the
	following properties: for $\Phi_{12} \in \{S_{12}, F_{12}\}$, we have that the diagram
	\[\begin{tikzcd}
		{\Delta(A, \nu_{M_1})} \arrow{r}{\Phi_{12}}\arrow{rd}{{p_{[1,d-1]} }}&{\Delta(A, \nu_{M_2})} \arrow{d}{{p_{[1,d-1]} }}\\
		 &{\Delta(\nu_M,A)}
	\end{tikzcd}\]
	commutes, and $\Phi_{12}$ preserves the Euclidean lengths of the fibers of $p_{[1,d-1]}$.
\end{theorem}

This theorem states that moving between two adjacent prime maximal cones in a tropical variety, corresponds to a volume preserving piece-wise linear map between their Newton-Okounkov bodies. Such maps are called mutations of polytopes, and they appear in many areas of combinatorics. The importance of this theorem is that one can travel from one prime cone to the next by the use of these mutations. This means it is possible to compute new-prime cones in the tropicalization starting from a known prime cone, by guessing the right mutations. In practice, this should be combined with knowledge about the form these mutations take on a certain tropical variety such as for $\Gr(2,n)$, which is done in \cite{WallCrossing} section 5.

\subsection{Wall-crossing maps}

The second part of \Cref{thm : wall-crossing}, pertaining to the wall crossing maps $S_{12}$ and $F_{12}$ is not very concrete.\\ 
However, these maps can be constructed explicitly, assuming knowledge of the inequality descriptions of the Newton-Okounkov bodies $\Delta(A,\nu_{M_1})$ and $\Delta(A,\nu_{M_2})$.

The authors call these maps the `geometric wall-crossing maps', as they also construct an `algebraic wall-crossing map' between the value semigroups $S(A,\nu_{M_1})$ and $S(A,\nu_{M_2})$. This algebraic wall-crossing map is not generally a nice map however, for example only in distinct cases is it a semigroup homomorphism or a restriction of one of the geometric wall-crossing maps. Therefore, for our purposes, it will not be worth discussing further.\\

For ease of notation, we will denote $S_i := S(A,\nu_{M_i})$ and $\Delta_i := \Delta(A,\nu_{M_i})$ in the rest of this section. Consider $\Delta_1$ and $\Delta_2$ again as subsets of $\{1\}\times\R^{d-1} \subset \R^d$.\\
It is clear that for $i \in \{1,2\}$ it is possible to find two piece-wise linear maps $\psi_i$ and $\phi_i$ from $\Delta:=\Delta(A,\nu_M)$ to $\R$ such that we have \[\Delta_i = \{(x,h)\in \R^{d-1}\times\R|x \in \Delta, \phi_i(x)\le h \le \psi_i(x)\}\]
Now let $\kappa$ be the global constant relating the Euclidean length of the fibers: \[\kappa(\psi_1(x) -\phi_1(x)) =\psi_2(x) -\phi_2(x) \]
We now define $S_{12}$ and $F_{12}$ as follows:
\[S_{12}(x,h) = (x, \kappa(h-\phi_1(x))+\phi_2(x))\]
\[F_{12}(x,h) = (x, -\kappa(h-\phi_1(x))+\psi_2(x))\]

It is clear these maps are piece-wise linear and make the required diagram commute, and an easy computation shows that they are bijective maps between $\Delta_1$ and $\Delta_2$.%

  \graphicspath{{\chaptersdir/Wall-crossing_for_non-prime_cones/image/}}%
  \chapter{Wall-crossing for non-prime cones}\label{ch:Wall-crossing for non-prime cones}

In this chapter we propose a method for wall crossing from a prime cone to an adjacent non-prime cone. The idea is to use an algorithm from \cite{Algorithm} to compute a new embedding of the variety in order to change the tropicalization in such a way that the non-prime cone becomes prime. The original prime cone will stay prime in this new embedding 
, and therefore if both cones are still adjacent in the new embedding, one can then apply the wall-crossing result from \cite{WallCrossing} and project back down to the original embedding, effectively resulting in a wall-crossing map for the original cones.

We then apply this method to an example in $\trop (\Gr(3,6))$, and show that it works in this case. For this we start from computations done in \cite{MatchingFields}, which sparked the idea of non-prime wall-crossing being possible in this example.

\section{Algorithm for computing a new embedding}
We start by detailing the algorithm from \cite{Algorithm}. In order to apply this algorithm to a maximal non-prime cone $C \subset \trop(I)$, we need the assumption that $C$ has multiplicity one. The formal definition of multiplicity of a cone can be found in \cite{IntroToTropGeom}, but we will define it using \cite[Lemma 1]{Algorithm} . First recall the following notions:
\begin{defn}[Primary ideal]
    An ideal $Q$ is called primary if the following holds: if $fg \in Q$ then $f\in Q$ or $g^n\in Q$ for some $n>0$. 
\end{defn}
It clearly follows from this definition that the radical ideal of a primary ideal is a prime ideal.

Now for any ideal $I$, it is possible to write $I=\bigcap_iQ_i$ where the $Q_i$ are a finite amount of primary ideals, with radical ideals $P_i = \sqrt{Q_i}$. Such a primary decomposition is not always unique, but the set $\{P_i\}$ is completely determined by $I$, and is called the set of associated primes of $I$.

We can now state what it means for a cone to have multiplicity one.

\begin{defn}[Maximal cone of multiplicity 1]
    If $I$ is a homogeneous ideal, we say a maximal cone $C\subset \trop(I)$ has multiplicity one if $\initial_C(I)$ has a unique toric associated prime ideal.
\end{defn}

Clearly, any prime maximal cone $C\subset \trop(I)$ has multiplicity equal to one, as the only associated prime of $\initial_C(I)$ is $\initial_C(I)$ itself. For a non-prime maximal cone however, this is not always the case.

The reason we work with the case of non-prime cones with multiplicity one is that, intuitively, a non-prime cone with multiplicity one is 'almost prime'. By this we mean that having multiplicity one is the closest a non-prime cone $C$ can get to being prime, as even though $\initial_C(I)$ is not toric, we can associate a unique toric ideal to $C$ by our definition of multiplicity one. In fact, in \cite{Algorithm}, the authors show that the assumption of multiplicity one  allows one to compute a new embedding of $V(I)$ such that, very informally said, the cone $C$ becomes prime.

The assumption that the non-prime cone we are working with has multiplicity one is not so restrictive. For example, in \cite{Algorithm}, multiplicities of maximal cones in the tropical varieties of several flag varieties are computed, and are all equal to one. 
\vspace{5pt}

We now show how to compute a new embedding. Note that the authors of \cite{Algorithm} did not develop this method with the wall-crossing theorem in mind, and did not concern themselves with what happens to adjacent cones in the new embedding. We will show later that at least in the case of $\Gr(3,6)$, things work out as we desire.

Let $I$ be an ideal in a polynomial ring $\C[x_1, ..., x_n]$, and let $C_1\subset \trop(I)$ be a maximal prime cone adjacent to a non-prime maximal cone $C_2\subset \trop(I)$ of multiplicity one. Here we consider $I$ in the Laurent polynomial ring $\C[x_1^{\pm1}, ..., x_n^{\pm1}]$ in order to compute $\trop(I)$. The original algorithm goes as follows:

First, compute a primary decomposition of $\initial_{C_2}(I)$, and compute its unique toric associated prime ideal $I_T$. This can be done by use of \cite[Lemma 2]{Algorithm}.\\
Now fix a minimal binomial generating set $G$ of $I_T$, and find the binomials that are in $G$ but not in $\initial_{C_2}(I)$, denote them by $L_{C_2} = \{f_1, ..., f_s\}$. We construct a new ideal $I' = I + \langle f_1-y_1, ..., f_s-y_s\rangle$ inside the polynomial ring $ \C[x_1,...,x_n,y_1,...,y_s]$. Since $ \C[x_1,...,x_n]/I \cong \C[x_1,...,x_n,y_1,...,y_s]/I'$, this is an embedding of the same variety, but in higher dimensions.\\
As such, any cone in the tropicalization of $I'$ can be projected onto the first $n$ coordinates to obtain a cone which is a subset of $\trop(I)$. If there is a prime maximal cone $C_2'$ in the tropicalization of $I'$ such that the projection of $C_2'$ onto the first $n$ coordinates  has a non-empty intersection with 
$C_2$, we stop as this is exactly what we want. If not, we repeat the process.
Note however that this algorithm does not always terminate as shown in \cite{ilten2020khovanskiifinite}.

\begin{algorithm}\label{algorithm:NewEmbedding1}
    \SetAlgoLined
    \SetKwInOut{Data}{Input}
    \SetKwInOut{Result}{Output}
    \Data{$A\cong \C[x_1, ..., x_n]/I$, and a non-prime maximal cone $C \subset \trop(I)$}
    \Result{An ideal $I'\subset R = \C[x_1,...,x_n, y_1,...,y_k]$ such that $A \cong R/I'$ and a prime maximal cone $C' \subset \trop(I')$ that projects onto $C$}
    Compute $I_T$: the unique toric associated prime of $\initial_C(I)$\;
    Fix $G$ = minimal binomial generating set of $I_T$\;
    Compute $L_C=\{f_1,...,f_s\}$ the list of binomials in $G$ but not in $\initial_C(I)$\;
    Set $I_1 = I+ \langle y_1-f_1, ..., y_s-f_s\rangle$ and $A_1 = \C[x_1,...,x_n, y_1,...,y_s]/I_1$ \;
    Compute $\trop(I_1)$\;
    \ForEach{Prime maximal cone $C_1$ in $\trop(I_1)$}{
        \eIf{$C_1 \text{ projects onto } C$}{
            return $I' = I_1, A' = A_1$ and $C' = C_1$\;
        } 
        {
            Apply algorithm 1 to  $A_1,I_1, C_1$ \;
        }
    }
    \caption{Computing a new embedding}
\end{algorithm}
\remark{In this context, when we say a cone projects onto another cone, we mean the projection of the first cone has a non-empty intersection with the relative interior of the second cone. This means there exists a weight vector in the first cone that, when projected onto the original embedding, yields the initial ideal of the second cone.}
It is possible to translate this algorithm into the language of Khovanskii bases. In that context, we are extending the valuation $\nu_C$ on $A$ to a valuation $\nu_{C_1}$ on $A_1$, such that the images of $x_1, ..., x_n, y_1, ..., y_n$ in $A_1$ form a Khovanskii basis for $(A_1, \nu_{C_1})$,see \cite[Example 7]{Algorithm}. \\
For our purposes, we also want a prime maximal cone $C_1'$ projecting onto $C_1$, such that $C_1'$ and $C_2'$ are adjacent. Therefore if we find a prime maximal cone projecting into the non-prime cone, we should check if there is an adjacent prime maximal cone projecting into the original prime maximal cone. Our adjusted algorithm would therefore look like this:

\begin{algorithm}\label{algorithm:NewEmbedding2}
    \SetAlgoLined
    \SetKwInOut{Data}{Input}
    \SetKwInOut{Result}{Output}
    \Data{$A\cong \C[x_1, ..., x_n]/I$,  a non-prime maximal cone $C_1 \subset \trop(I)$ and an adjacent prime maximal cone $C_2\subset \trop(I)$}
    \Result{An ideal $I'\subset R = \C[x_1,...,x_n, y_1,...,y_k]$ such that $A \cong R/I'$, a prime maximal cone $C_1' \subset \trop(I')$ that projects onto $C_1$ and an adjacent prime maximal cone $C_2'$ that projects onto $C_2$}
    Compute $I_T$: the unique toric associated prime of $\initial_{C_1}(I)$\;
    Fix $G$ = minimal binomial generating set of $I_T$\;
    Compute $L_{C_1}=\{f_1,...,f_s\}$ the list of binomials in $G$ but not in $\initial_C(I)$\;
    Set $I_1 = I+ \langle y_1-f_1, ..., y_s-f_s\rangle$ and $A_1 = \C[x_1,...,x_n, y_1,...,y_s]/I_1$ \;
    Compute $\trop(I_1)$\;
    \ForEach{Prime maximal cone $C_1'$ in $\trop(I_1)$}{
        \eIf{$C_1' \text{ projects into } C_1$}{\ForEach{Prime maximal cone $C_2'$ adjacent to $C_1'$}{\If{$C_2' \text{ projects onto } C_2$}{
            return $I' = I_1, A' = A_1$, $C_1'$ and $C_2'$\;
        } }}
        {
            Apply algorithm 1 to  $A_1,I_1, C_1$ \;
        }
    }
    \caption{Adjustment of algorithm 1 for adjacent cones}
\end{algorithm}
It is important to stress that it is unsure if we can always find such an adjacent cone $C_2'$. However, if we do find such cones $C_1'$ and $C_2'$, we can apply the wall-crossing theorem to them to obtain a mutation between their associated Newton-Okounkov bodies. This would allow us to wall-cross from $C_1$ to $C_2$ after projecting back down. 

In the next section, we show that this approach works for $\Gr(3,6)$, meaning we can find such adjacent cones in the new embedding. In fact, what we find is an even better result than expected.

\section{Computations for Gr(3,6)}

In this section we detail our computations for applying the above procedure to certain cones in $\Gr(3,6)$. The computations are done in \textit{Macaulay2} and the code can be found at \url{https://github.com/CasProost/NObodiesToricDegens}.

The cones we look at come from computations done in \cite{MatchingFields} in the context of matching fields. Matching fields are combinatorial objects that have been successfully used to define toric degenerations of Grassmannian, flag varieties and Schubert varieties \cite{ToricdegensofGrassmannians, clarke2020toric, OllieFatemeh3, OllieFatemeh2, MatchingFields}.
In \cite[Example 4]{MatchingFields} the authors obtain, through the use of these matching fields, a mutation of a Newton-Okounkov polytope associated to $\Gr(3,6)$ to another Newton-Okounkov polytope of $\Gr(3,6)$. Upon closer inspection the authors found that this mutation is a composition of three mutations, where after the first  two mutations in this composition, a non-lattice polytope is obtained. This non-lattice polytope can not come from a prime cone, and hence this seemed to be wall-crossing through a non-prime cone.

Below we translate this example to our language of the tropical Grassmannian $\trop(\Gr(3,6))$, apply the algorithm to compute a new embedding and show that the necessary adjacent prime maximal cones exist in the new embedding. This was all done by making use of the computer algebra software Macaulay2 \cite{M2}. 

\subsection{Finding the cones}

We follow the labeling from the original source of the example. The starting matching field is called $\sigma$, with matching field polytope $\Delta_\sigma$, and it is mutated to a matching field polytope $\Delta_\tau$ corresponding to a matching field $\tau$. We will call the polytope obtained after the first of the three mutations $\Delta_{\rho}$ and the non-lattice polytope obtained after the first two mutations ($Q$ in the original paper) $\Delta $. The intermediate polytope $\Delta_\rho$ also corresponds to a matching  field $\rho$ and a prime maximal cone. As such the first mutation, from $\Delta_\sigma$ to $\Delta_\rho$ corresponds to a wall-crossing between adjacent prime cones.

The first step is to compute $\trop(\Gr(3,6))$ and find the maximal cones corresponding to $\Delta_\sigma, \Delta_\tau, \Delta_{\rho}$ and $\Delta $. Using the Macaulay2 package \cite{MatchingFieldsSource} it is possible to obtain a weight vector in $\trop(\Gr(3,6))$ corresponding to a matching field, which allows us to check which maximal prime cones correspond to $\Delta_\sigma, \Delta_\tau$ and $\Delta_{\rho}$. We call these cones $C_\sigma$, $C_\tau$ and $C_{\rho}$ respectively. We include the rays of these cones as the column vectors of the following matrix: 
\[\left(\!\begin{array}{ccccccc}
11&21&1&-2&1&3&3\\
1&-4&1&3&11&-2&3\\
11&-4&1&-2&1&3&3\\
1&-4&6&3&11&-2&3\\
-9&-4&-4&-2&-9&-2&-2\\
1&-4&-4&3&11&3&-2\\
-9&-4&1&-2&-9&-2&-2\\
-9&1&-4&-2&-9&-2&-2\\
11&1&1&3&1&3&-2\\
-9&1&1&-2&-9&-2&-2\\
-9&-4&-4&-2&-9&-2&-2\\
1&-4&-4&3&11&3&-2\\
-9&-4&1&-2&-9&-2&-2\\
-9&1&-4&-2&-9&-2&-2\\
11&1&1&3&1&3&-2\\
-9&1&1&-2&-9&-2&-2\\
11&1&21&3&1&-2&3\\
1&1&-4&-2&11&3&3\\
11&1&-4&3&1&-2&3\\
1&6&-4&-2&11&3&3
\end{array}\!\right)
 \]

We will denote $r_i$ for the $i$-th column of this matrix. We now give the rays which (together with the lineality space) generate each cone: \[C_\sigma: r_1, r_2, r_3, r_7\]
\[C_\tau: r_1, r_2, r_6, r_7\]
\[C_\rho: r_1, r_3, r_4, r_7\]

Since the non-lattice polytope does not arise from a matching field, we can not immediately compute a weight vector corresponding to this polytope to find the corresponding non-prime cone. Instead, we look for non-prime cones adjacent to $C_{\rho}$ and $C_\tau$. 

Unfortunately, there does not seem to be a single non-prime cone that is adjacent to both $C_{\rho}$ and $C_\tau$. It seems like the example from \cite{MatchingFields} does not simply correspond to going from a prime cone, through a non-prime cone, to another prime cone. What we did find are two adjacent maximal non-prime cones, one adjacent to $C_\rho$ and one adjacent to $C_\tau$. We will denote the non-prime cone adjacent to $C_\rho$ by $C_1$ and the non-prime cone adjacent to $C_\tau$ by $C_2$. Their rays are:

\[C_1: r_1, r_4, r_5, r_7\]
\[C_1: r_1, r_5, r_6, r_7\]

A priori, it might seem like we have pulled these non-prime cones out of thin air, however, our computations show that we have good reason to believe these cones correspond, together, to the non-lattice polytope $\Delta$.

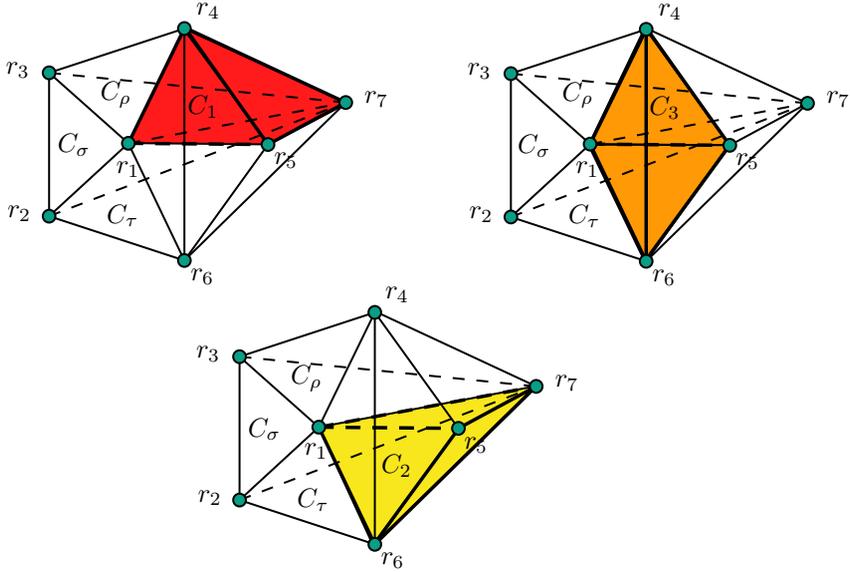
\begin{figure}
    \begin{center}
        \tikzset{every picture/.style={line width=0.75pt}} 

\begin{tikzpicture}[x=0.6pt,y=0.6pt,yscale=-1,xscale=1]

\draw  [draw opacity=0][fill={rgb, 255:red, 248; green, 231; blue, 28 }  ,fill opacity=0.5 ] (239.5,293) -- (375.89,267.42) -- (273.74,366.93) -- cycle ;
\draw [line width=1.3]    (273.74,366.93) -- (327.32,294.29) ;
\draw  [draw opacity=0][fill={rgb, 255:red, 255; green, 153; blue, 5 }  ,fill opacity=0.5 ] (444.92,43.3) -- (498.5,116.75) -- (409.91,115.92) -- cycle ;
\draw  [draw opacity=0][fill={rgb, 255:red, 255; green, 153; blue, 5 }  ,fill opacity=0.5 ] (444.91,189.39) -- (409.33,116.31) -- (497.91,116.67) -- cycle ;
\draw  [draw opacity=0][fill={rgb, 255:red, 255; green, 27; blue, 27 }  ,fill opacity=0.5 ] (154.81,42.73) -- (208.4,116.18) -- (119.81,115.35) -- cycle ;
\draw  [draw opacity=0][fill={rgb, 255:red, 255; green, 27; blue, 27 }  ,fill opacity=0.5 ] (207.99,116.12) -- (153.68,43.22) -- (257.49,90.33) -- cycle ;
\draw [color={rgb, 255:red, 0; green, 0; blue, 0 }  ,draw opacity=1 ][line width=1.3]    (207.5,116.25) -- (154.81,42.73) ;
\draw [color={rgb, 255:red, 0; green, 0; blue, 0 }  ,draw opacity=0.7 ]   (154.81,43.73) -- (154.81,188.77) ;
\draw  [dash pattern={on 4.5pt off 4.5pt}]  (119.5,115.5) -- (256.5,90) ;
\draw [color={rgb, 255:red, 0; green, 0; blue, 0 }  ,draw opacity=0.7 ] [dash pattern={on 4.5pt off 4.5pt}]  (69.56,161.07) -- (256.5,90) ;
\draw [color={rgb, 255:red, 0; green, 0; blue, 0 }  ,draw opacity=0.7 ] [dash pattern={on 4.5pt off 4.5pt}]  (69.56,71.43) -- (256.5,90) ;
\draw [line width=1.3]    (207.5,116.25) -- (256.5,90) ;
\draw [color={rgb, 255:red, 0; green, 0; blue, 0 }  ,draw opacity=0.7 ]   (154.81,188.77) -- (256.5,90) ;
\draw [line width=1.3]    (154.5,42.89) -- (256.5,90) ;
\draw  [color={rgb, 255:red, 0; green, 0; blue, 0 }  ,draw opacity=0.7 ] (207.5,116.25) -- (154.81,188.77) -- (69.56,161.07) -- (69.56,71.43) -- (154.81,43.73) -- cycle ;
\draw [color={rgb, 255:red, 0; green, 0; blue, 0 }  ,draw opacity=0.7 ]   (69.56,71.43) -- (119.5,115.5) ;
\draw [color={rgb, 255:red, 0; green, 0; blue, 0 }  ,draw opacity=0.7 ]   (119.5,115.5) -- (154.81,188.77) ;
\draw [line width=1.3]  [dash pattern={on 5.63pt off 4.5pt}]  (119.5,115.5) -- (207.5,116.25) ;
\draw [color={rgb, 255:red, 0; green, 0; blue, 0 }  ,draw opacity=1 ][line width=1.3]    (119.5,115.5) -- (154.81,43.73) ;
\draw [color={rgb, 255:red, 0; green, 0; blue, 0 }  ,draw opacity=0.7 ]   (119.5,115.5) -- (69.56,161.07) ;
\draw  [fill={rgb, 255:red, 9; green, 158; blue, 134 }  ,fill opacity=1 ] (115.5,115.5) .. controls (115.5,113.29) and (117.29,111.5) .. (119.5,111.5) .. controls (121.71,111.5) and (123.5,113.29) .. (123.5,115.5) .. controls (123.5,117.71) and (121.71,119.5) .. (119.5,119.5) .. controls (117.29,119.5) and (115.5,117.71) .. (115.5,115.5) -- cycle ;
\draw  [fill={rgb, 255:red, 9; green, 158; blue, 134 }  ,fill opacity=1 ] (150.81,43.73) .. controls (150.81,41.52) and (152.6,39.73) .. (154.81,39.73) .. controls (157.02,39.73) and (158.81,41.52) .. (158.81,43.73) .. controls (158.81,45.94) and (157.02,47.73) .. (154.81,47.73) .. controls (152.6,47.73) and (150.81,45.94) .. (150.81,43.73) -- cycle ;
\draw  [fill={rgb, 255:red, 9; green, 158; blue, 134 }  ,fill opacity=1 ] (203.5,116.25) .. controls (203.5,114.04) and (205.29,112.25) .. (207.5,112.25) .. controls (209.71,112.25) and (211.5,114.04) .. (211.5,116.25) .. controls (211.5,118.46) and (209.71,120.25) .. (207.5,120.25) .. controls (205.29,120.25) and (203.5,118.46) .. (203.5,116.25) -- cycle ;
\draw  [fill={rgb, 255:red, 9; green, 158; blue, 134 }  ,fill opacity=1 ] (150.81,188.77) .. controls (150.81,186.56) and (152.6,184.77) .. (154.81,184.77) .. controls (157.02,184.77) and (158.81,186.56) .. (158.81,188.77) .. controls (158.81,190.98) and (157.02,192.77) .. (154.81,192.77) .. controls (152.6,192.77) and (150.81,190.98) .. (150.81,188.77) -- cycle ;
\draw  [fill={rgb, 255:red, 9; green, 158; blue, 134 }  ,fill opacity=1 ] (65.56,161.07) .. controls (65.56,158.86) and (67.35,157.07) .. (69.56,157.07) .. controls (71.77,157.07) and (73.56,158.86) .. (73.56,161.07) .. controls (73.56,163.28) and (71.77,165.07) .. (69.56,165.07) .. controls (67.35,165.07) and (65.56,163.28) .. (65.56,161.07) -- cycle ;
\draw  [fill={rgb, 255:red, 9; green, 158; blue, 134 }  ,fill opacity=1 ] (65.56,71.43) .. controls (65.56,69.22) and (67.35,67.43) .. (69.56,67.43) .. controls (71.77,67.43) and (73.56,69.22) .. (73.56,71.43) .. controls (73.56,73.64) and (71.77,75.43) .. (69.56,75.43) .. controls (67.35,75.43) and (65.56,73.64) .. (65.56,71.43) -- cycle ;
\draw  [fill={rgb, 255:red, 9; green, 158; blue, 134 }  ,fill opacity=1 ] (252.5,90) .. controls (252.5,87.79) and (254.29,86) .. (256.5,86) .. controls (258.71,86) and (260.5,87.79) .. (260.5,90) .. controls (260.5,92.21) and (258.71,94) .. (256.5,94) .. controls (254.29,94) and (252.5,92.21) .. (252.5,90) -- cycle ;
\draw [line width=1.3]    (444.91,189.39) -- (498.5,116.75) ;
\draw [color={rgb, 255:red, 0; green, 0; blue, 0 }  ,draw opacity=1 ][line width=1.3]    (498.5,116.75) -- (445.81,43.23) ;
\draw [line width=1.3]    (445.81,44.23) -- (445.81,189.27) ;
\draw  [dash pattern={on 4.5pt off 4.5pt}]  (410.5,116) -- (547.5,90.5) ;
\draw [color={rgb, 255:red, 0; green, 0; blue, 0 }  ,draw opacity=0.7 ] [dash pattern={on 4.5pt off 4.5pt}]  (360.56,161.57) -- (547.5,90.5) ;
\draw [color={rgb, 255:red, 0; green, 0; blue, 0 }  ,draw opacity=0.7 ] [dash pattern={on 4.5pt off 4.5pt}]  (360.56,71.93) -- (547.5,90.5) ;
\draw [color={rgb, 255:red, 0; green, 0; blue, 0 }  ,draw opacity=0.7 ][line width=0.75]    (498.5,116.75) -- (547.5,90.5) ;
\draw [color={rgb, 255:red, 0; green, 0; blue, 0 }  ,draw opacity=0.7 ]   (445.81,189.27) -- (547.5,90.5) ;
\draw [line width=0.75]    (445.5,43.39) -- (547.5,90.5) ;
\draw  [color={rgb, 255:red, 0; green, 0; blue, 0 }  ,draw opacity=0.7 ] (498.5,116.75) -- (445.81,189.27) -- (360.56,161.57) -- (360.56,71.93) -- (445.81,44.23) -- cycle ;
\draw [color={rgb, 255:red, 0; green, 0; blue, 0 }  ,draw opacity=0.7 ]   (360.56,71.93) -- (410.5,116) ;
\draw [line width=1.3]    (410.5,116) -- (445.81,189.27) ;
\draw [line width=1.3]  [dash pattern={on 5.63pt off 4.5pt}]  (410.5,116) -- (498.5,116.75) ;
\draw [color={rgb, 255:red, 0; green, 0; blue, 0 }  ,draw opacity=1 ][line width=1.3]    (410.5,116) -- (445.81,44.23) ;
\draw [color={rgb, 255:red, 0; green, 0; blue, 0 }  ,draw opacity=0.7 ]   (410.5,116) -- (360.56,161.57) ;
\draw  [fill={rgb, 255:red, 9; green, 158; blue, 134 }  ,fill opacity=1 ] (406.5,116) .. controls (406.5,113.79) and (408.29,112) .. (410.5,112) .. controls (412.71,112) and (414.5,113.79) .. (414.5,116) .. controls (414.5,118.21) and (412.71,120) .. (410.5,120) .. controls (408.29,120) and (406.5,118.21) .. (406.5,116) -- cycle ;
\draw  [fill={rgb, 255:red, 9; green, 158; blue, 134 }  ,fill opacity=1 ] (441.81,44.23) .. controls (441.81,42.02) and (443.6,40.23) .. (445.81,40.23) .. controls (448.02,40.23) and (449.81,42.02) .. (449.81,44.23) .. controls (449.81,46.44) and (448.02,48.23) .. (445.81,48.23) .. controls (443.6,48.23) and (441.81,46.44) .. (441.81,44.23) -- cycle ;
\draw  [fill={rgb, 255:red, 9; green, 158; blue, 134 }  ,fill opacity=1 ] (494.5,116.75) .. controls (494.5,114.54) and (496.29,112.75) .. (498.5,112.75) .. controls (500.71,112.75) and (502.5,114.54) .. (502.5,116.75) .. controls (502.5,118.96) and (500.71,120.75) .. (498.5,120.75) .. controls (496.29,120.75) and (494.5,118.96) .. (494.5,116.75) -- cycle ;
\draw  [fill={rgb, 255:red, 9; green, 158; blue, 134 }  ,fill opacity=1 ] (441.81,189.27) .. controls (441.81,187.06) and (443.6,185.27) .. (445.81,185.27) .. controls (448.02,185.27) and (449.81,187.06) .. (449.81,189.27) .. controls (449.81,191.48) and (448.02,193.27) .. (445.81,193.27) .. controls (443.6,193.27) and (441.81,191.48) .. (441.81,189.27) -- cycle ;
\draw  [fill={rgb, 255:red, 9; green, 158; blue, 134 }  ,fill opacity=1 ] (356.56,161.57) .. controls (356.56,159.36) and (358.35,157.57) .. (360.56,157.57) .. controls (362.77,157.57) and (364.56,159.36) .. (364.56,161.57) .. controls (364.56,163.78) and (362.77,165.57) .. (360.56,165.57) .. controls (358.35,165.57) and (356.56,163.78) .. (356.56,161.57) -- cycle ;
\draw  [fill={rgb, 255:red, 9; green, 158; blue, 134 }  ,fill opacity=1 ] (356.56,71.93) .. controls (356.56,69.72) and (358.35,67.93) .. (360.56,67.93) .. controls (362.77,67.93) and (364.56,69.72) .. (364.56,71.93) .. controls (364.56,74.14) and (362.77,75.93) .. (360.56,75.93) .. controls (358.35,75.93) and (356.56,74.14) .. (356.56,71.93) -- cycle ;
\draw  [fill={rgb, 255:red, 9; green, 158; blue, 134 }  ,fill opacity=1 ] (543.5,90.5) .. controls (543.5,88.29) and (545.29,86.5) .. (547.5,86.5) .. controls (549.71,86.5) and (551.5,88.29) .. (551.5,90.5) .. controls (551.5,92.71) and (549.71,94.5) .. (547.5,94.5) .. controls (545.29,94.5) and (543.5,92.71) .. (543.5,90.5) -- cycle ;
\draw [color={rgb, 255:red, 0; green, 0; blue, 0 }  ,draw opacity=0.7 ]   (274.81,221.23) -- (274.81,366.27) ;
\draw [line width=1.3]  [dash pattern={on 5.63pt off 4.5pt}]  (239.5,293) -- (376.5,267.5) ;
\draw [color={rgb, 255:red, 0; green, 0; blue, 0 }  ,draw opacity=0.7 ] [dash pattern={on 4.5pt off 4.5pt}]  (189.56,338.57) -- (376.5,267.5) ;
\draw [color={rgb, 255:red, 0; green, 0; blue, 0 }  ,draw opacity=0.7 ] [dash pattern={on 4.5pt off 4.5pt}]  (189.56,248.93) -- (376.5,267.5) ;
\draw [line width=1.3]    (327.5,293.75) -- (376.5,267.5) ;
\draw [line width=1.3]    (274.81,366.27) -- (376.5,267.5) ;
\draw [line width=0.75]    (274.5,220.39) -- (376.5,267.5) ;
\draw  [color={rgb, 255:red, 0; green, 0; blue, 0 }  ,draw opacity=0.7 ] (327.5,293.75) -- (274.81,366.27) -- (189.56,338.57) -- (189.56,248.93) -- (274.81,221.23) -- cycle ;
\draw [color={rgb, 255:red, 0; green, 0; blue, 0 }  ,draw opacity=0.7 ]   (189.56,248.93) -- (239.5,293) ;
\draw [line width=1.3]    (239.5,293) -- (274.81,366.27) ;
\draw [line width=1.3]  [dash pattern={on 5.63pt off 4.5pt}]  (239.5,293) -- (327.5,293.75) ;
\draw [color={rgb, 255:red, 0; green, 0; blue, 0 }  ,draw opacity=1 ][line width=0.75]    (239.5,293) -- (274.81,221.23) ;
\draw [color={rgb, 255:red, 0; green, 0; blue, 0 }  ,draw opacity=0.7 ]   (239.5,293) -- (189.56,338.57) ;
\draw  [fill={rgb, 255:red, 9; green, 158; blue, 134 }  ,fill opacity=1 ] (235.5,293) .. controls (235.5,290.79) and (237.29,289) .. (239.5,289) .. controls (241.71,289) and (243.5,290.79) .. (243.5,293) .. controls (243.5,295.21) and (241.71,297) .. (239.5,297) .. controls (237.29,297) and (235.5,295.21) .. (235.5,293) -- cycle ;
\draw  [fill={rgb, 255:red, 9; green, 158; blue, 134 }  ,fill opacity=1 ] (270.81,221.23) .. controls (270.81,219.02) and (272.6,217.23) .. (274.81,217.23) .. controls (277.02,217.23) and (278.81,219.02) .. (278.81,221.23) .. controls (278.81,223.44) and (277.02,225.23) .. (274.81,225.23) .. controls (272.6,225.23) and (270.81,223.44) .. (270.81,221.23) -- cycle ;
\draw  [fill={rgb, 255:red, 9; green, 158; blue, 134 }  ,fill opacity=1 ] (323.5,293.75) .. controls (323.5,291.54) and (325.29,289.75) .. (327.5,289.75) .. controls (329.71,289.75) and (331.5,291.54) .. (331.5,293.75) .. controls (331.5,295.96) and (329.71,297.75) .. (327.5,297.75) .. controls (325.29,297.75) and (323.5,295.96) .. (323.5,293.75) -- cycle ;
\draw  [fill={rgb, 255:red, 9; green, 158; blue, 134 }  ,fill opacity=1 ] (270.81,366.27) .. controls (270.81,364.06) and (272.6,362.27) .. (274.81,362.27) .. controls (277.02,362.27) and (278.81,364.06) .. (278.81,366.27) .. controls (278.81,368.48) and (277.02,370.27) .. (274.81,370.27) .. controls (272.6,370.27) and (270.81,368.48) .. (270.81,366.27) -- cycle ;
\draw  [fill={rgb, 255:red, 9; green, 158; blue, 134 }  ,fill opacity=1 ] (185.56,338.57) .. controls (185.56,336.36) and (187.35,334.57) .. (189.56,334.57) .. controls (191.77,334.57) and (193.56,336.36) .. (193.56,338.57) .. controls (193.56,340.78) and (191.77,342.57) .. (189.56,342.57) .. controls (187.35,342.57) and (185.56,340.78) .. (185.56,338.57) -- cycle ;
\draw  [fill={rgb, 255:red, 9; green, 158; blue, 134 }  ,fill opacity=1 ] (185.56,248.93) .. controls (185.56,246.72) and (187.35,244.93) .. (189.56,244.93) .. controls (191.77,244.93) and (193.56,246.72) .. (193.56,248.93) .. controls (193.56,251.14) and (191.77,252.93) .. (189.56,252.93) .. controls (187.35,252.93) and (185.56,251.14) .. (185.56,248.93) -- cycle ;
\draw  [fill={rgb, 255:red, 9; green, 158; blue, 134 }  ,fill opacity=1 ] (372.5,267.5) .. controls (372.5,265.29) and (374.29,263.5) .. (376.5,263.5) .. controls (378.71,263.5) and (380.5,265.29) .. (380.5,267.5) .. controls (380.5,269.71) and (378.71,271.5) .. (376.5,271.5) .. controls (374.29,271.5) and (372.5,269.71) .. (372.5,267.5) -- cycle ;

\draw (73,107.9) node [anchor=north west][inner sep=0.75pt]  [color={rgb, 255:red, 0; green, 0; blue, 0 }  ,opacity=0.7 ]  {$C_{\sigma }$};
\draw (104,151.9) node [anchor=north west][inner sep=0.75pt]  [color={rgb, 255:red, 0; green, 0; blue, 0 }  ,opacity=0.7 ]  {$C_{\tau }$};
\draw (100,74.9) node [anchor=north west][inner sep=0.75pt]  [color={rgb, 255:red, 0; green, 0; blue, 0 }  ,opacity=0.7 ]  {$C_{\rho }$};
\draw (110,125) node [anchor=north west][inner sep=0.75pt]    {$r_{1}$};
\draw (41.5,152.9) node [anchor=north west][inner sep=0.75pt]  [color={rgb, 255:red, 0; green, 0; blue, 0 }  ,opacity=0.7 ]  {$r_{2}$};
\draw (40.5,62.9) node [anchor=north west][inner sep=0.75pt]  [color={rgb, 255:red, 0; green, 0; blue, 0 }  ,opacity=0.7 ]  {$r_{3}$};
\draw (160.5,25.9) node [anchor=north west][inner sep=0.75pt]    {$r_{4}$};
\draw (209.5,119.65) node [anchor=north west][inner sep=0.75pt]    {$r_{5}$};
\draw (156.81,192.17) node [anchor=north west][inner sep=0.75pt]  [color={rgb, 255:red, 0; green, 0; blue, 0 }  ,opacity=0.7 ]  {$r_{6}$};
\draw (266.5,80.9) node [anchor=north west][inner sep=0.75pt]    {$r_{7}$};
\draw (155,83.9) node [anchor=north west][inner sep=0.75pt]    {$C_{1}$};
\draw (363,108.4) node [anchor=north west][inner sep=0.75pt]  [color={rgb, 255:red, 0; green, 0; blue, 0 }  ,opacity=0.7 ]  {$C_{\sigma }$};
\draw (395,152.4) node [anchor=north west][inner sep=0.75pt]  [color={rgb, 255:red, 0; green, 0; blue, 0 }  ,opacity=0.7 ]  {$C_{\tau }$};
\draw (391,75.4) node [anchor=north west][inner sep=0.75pt]  [color={rgb, 255:red, 0; green, 0; blue, 0 }  ,opacity=0.7 ]  {$C_{\rho }$};
\draw (399.5,125) node [anchor=north west][inner sep=0.75pt]    {$r_{1}$};
\draw (332.5,153.4) node [anchor=north west][inner sep=0.75pt]  [color={rgb, 255:red, 0; green, 0; blue, 0 }  ,opacity=0.7 ]  {$r_{2}$};
\draw (331.5,63.4) node [anchor=north west][inner sep=0.75pt]  [color={rgb, 255:red, 0; green, 0; blue, 0 }  ,opacity=0.7 ]  {$r_{3}$};
\draw (451.5,26.4) node [anchor=north west][inner sep=0.75pt]    {$r_{4}$};
\draw (500.5,120.15) node [anchor=north west][inner sep=0.75pt]    {$r_{5}$};
\draw (447.81,192.67) node [anchor=north west][inner sep=0.75pt]    {$r_{6}$};
\draw (557.5,81.4) node [anchor=north west][inner sep=0.75pt]  [color={rgb, 255:red, 0; green, 0; blue, 0 }  ,opacity=0.7 ]  {$r_{7}$};
\draw (446,84.4) node [anchor=north west][inner sep=0.75pt]    {$C_{3}$};
\draw (193,285.4) node [anchor=north west][inner sep=0.75pt]  [color={rgb, 255:red, 0; green, 0; blue, 0 }  ,opacity=0.7 ]  {$C_{\sigma }$};
\draw (224,329.4) node [anchor=north west][inner sep=0.75pt]  [color={rgb, 255:red, 0; green, 0; blue, 0 }  ,opacity=0.7 ]  {$C_{\tau }$};
\draw (220,252.4) node [anchor=north west][inner sep=0.75pt]  [color={rgb, 255:red, 0; green, 0; blue, 0 }  ,opacity=0.7 ]  {$C_{\rho }$};
\draw (228.5,301) node [anchor=north west][inner sep=0.75pt]    {$r_{1}$};
\draw (161.5,330.4) node [anchor=north west][inner sep=0.75pt]  [color={rgb, 255:red, 0; green, 0; blue, 0 }  ,opacity=0.7 ]  {$r_{2}$};
\draw (160.5,240.4) node [anchor=north west][inner sep=0.75pt]  [color={rgb, 255:red, 0; green, 0; blue, 0 }  ,opacity=0.7 ]  {$r_{3}$};
\draw (279.5,203.4) node [anchor=north west][inner sep=0.75pt]  [color={rgb, 255:red, 0; green, 0; blue, 0 }  ,opacity=0.7 ]  {$r_{4}$};
\draw (329.32,297.69) node [anchor=north west][inner sep=0.75pt]    {$r_{5}$};
\draw (276.81,369.67) node [anchor=north west][inner sep=0.75pt]    {$r_{6}$};
\draw (386.5,258.4) node [anchor=north west][inner sep=0.75pt]    {$r_{7}$};
\draw (277,308.4) node [anchor=north west][inner sep=0.75pt]    {$C_{2}$};

\end{tikzpicture}
    \end{center}
    \caption{Representation of the relevant cones of $\trop \Gr(3,6)$, courtesy of F. Zaffalon.}
    \label{fig : conesGr36}
\end{figure}

The relevant cones are depicted in \cref{fig : conesGr36}. In this figure, one sees an arrangement of adjacent tetrahedra, each of which corresponds to a cone. This is a common way to depict or discuss higher dimensional tropical varieties. To understand how this works, recall that $\trop\Gr(3,6)$ is four dimensional when quotienting out the lineality space. One has to think of these tetrahedra as sitting inside four dimensional space, and the cone associated to one of them is the cone over this tetrahedron. For a polytope $P$, the cone over $P$ is obtained by placing the polytope some distance away from the origin, and then drawing rays from the origin to all points of $P$. Think of the triangle in $\R^3$ with vertices $(1,0,0), (0,1,0), (0,0,1)$, which yields the cone 
\[\text{cone} ((1,0,0), (0,1,0), (0,0,1))\] which is exactly the positive orthant. This is of course harder to visualize for the tetrahedra, as they sit in four dimensional space.\\
In \cref{fig : conesGr36} one sees that $C_\sigma$ and $C_\tau$ are actually adjacent. It might therefore seem like attempting to wall-cross through the non-prime cones is an unnecessary detour, and while there is some truth to that, this example can still serve as a proof of concept of Algorithm 2. Moreover, from  computations in \cite{MatchingFields}, it seems that finding a single mutation from $\Delta_\sigma$ to $\Delta_\tau$ is not possible.
Another interesting point is the cone $C_3$ appearing in this figure. It is given by the rays: \[C_3: r_1, r_4, r_5, r_6.\] This is a non-prime cone that is adjacent to both $C_1$ and $C_2$, that we originally did not consider. However when we compute the new embedding it will become clear that this cone is also relevant. Together, the tetrahedra of $C_1$, $C_2$ and $C_3$ form a bipyramid. In $\Gr(3,6)$, there are fifteen such bipyramids formed by three tetrahedra corresponding to non-prime cones \cite{TropicalGrassmannian}. We will discuss further in more detail.

\subsection{Computing a new embedding}

We start our computations by applying algorithm 2 to $C_\rho$ and $C_1$. First up is to compute the initial ideal of $I_{3,6}$ with respect to $C_1$. We get a non-prime  ideal, generated by binomials of degree 2, which we will denote by $I$. The list of generators can be recovered with Macaulay from the rays we gave above, and is quite unwieldy, and we therefore do not include it.

Now we can compute the new embedding. This is done by computing a primary decomposition of $I$, which yields its associated primes. There is a unique toric associated prime, so the assumption that $C_1$ is of multiplicity one is valid. We compute the set of missing binomials, which consists of the single binomial:
\[f = p_{126}p_{345}+p_{124}p_{356}\]

This allows us to compute the new embedding as \[I' = I_{3,6} =\langle Yh - f\rangle\] inside the polynomial ring $\C[p_I, Y,h]$. Here the variable $h$ is added for homogenization purposes, which we discuss below. Again using Macaulay2, we compute $\trop(I')$. 

\begin{theorem}
    The tropical variety $\trop(I')/V$ is an 4-dimensional fan with f-vector $\{1,78,692,1790,1337\}$, and $V$ is the lineality space of dimension $7$
\end{theorem}

\remark{The fact that the computation for this tropical variety finishes in Macaulay2 is not to be taken for granted. The runtime of this computation was anywhere between 45 minutes and 4 hours depending on the machine used. One usually would employ various known symmetries of the ideal to simplify the calculation for Macaulay2, but as we believe to be the first studying this expanded ideal, we could not rely on such known symmetries.}

We can compare this to the results about $\trop \Gr(3,6)$ from \cref{thm:tropGR(3;6)}, and we see that, when quotienting out the lineality space both tropical varieties have the same dimension. The new tropical variety however, contains 302 more maximal cones, which could indicate that some cones have been subdivided by the algorithm, or that some maximal cones project onto non-maximal cones.

Next up is checking if there are prime maximal cones projecting onto $C_1$, and if so, if there are adjacent prime maximal cones projecting onto $C_\rho$. This is where our results exceeded our expectations.

First of all we check for prime maximal cones projecting onto $C_1$. Our approach for this computation was as follows: for every cone $C$ in the new embedding, Macaulay gives us the list $L_C$ of rays that generate that cone. We project each of these rays down to the original embedding, by forgetting the last two coordinates (corresponding to $Y$ and $h$). The projection of $C$ is then the cone generated by the images of the rays in $L_C$.\\
Now we can check for every maximal cone $C$ in $\trop(I')$ if it projects onto $C_1$, and if it does, we compute the initial ideal $\initial_C(I')$. If this is a prime ideal, then $C$ is a prime maximal cone projecting onto $C_1$.

The result yields two such cones projecting onto $C_1$, let us call them $C_1'$ and $C_2'$. For each of these cones now, we can check if there are adjacent prime maximal cones projecting onto $C_\rho$. For this we first compute all prime maximal cones projecting onto $C_\rho$, and then check if they are adjacent to one of our two cones by computing the dimension of the intersection with either one of the cones. As we desire, there is in fact a maximal prime cone $C_\rho'\subset \trop(I')$ that projects onto $C_\rho$ and that is adjacent to one of the cones projecting onto $C_1$, which we shall assume to be $C_1'$.

This is a great result, as it means we can apply the wall-crossing theorem to cross between these two cones in the new embedding. However in order to pass through the non-prime cones in the original embedding, there would also have to be a prime cone in $\trop(I')$ that projects onto $C_2$, which would need to be adjacent to $C_1'$, and we would also like it to be adjacent to a prime cone $C_\tau'$ that projects onto $C_\tau$.

Remarkably, we already have this: $C_1'$ projects onto $C_2$, and there is a maximal prime cone $C_\tau'\subset \trop(I')$ that is adjacent to $C_1'$. Moreover, there is a maximal prime cone $C_\sigma' \subset \trop(I')$ that projects onto $C_\sigma$ and is adjacent to $C_\rho'$. This means that the entirety of the mutation from \cite{MatchingFields} is explained by this computation, only needing to compute one new embedding. Our expectations were that we would have to compute multiple new embeddings to cross to the multiple cones, but the situation worked out to be much nicer. 

\begin{figure}
    \begin{center}
        \tikzset{every picture/.style={line width=0.75pt}} 

\begin{tikzpicture}[x=0.6pt,y=0.6pt,yscale=-1,xscale=1]

\draw  [draw opacity=0][fill={rgb, 255:red, 80; green, 227; blue, 194 }  ,fill opacity=0.5 ] (555.82,80.73) -- (455.23,178.46) -- (454.25,34.01) -- cycle ;
\draw  [draw opacity=0][fill={rgb, 255:red, 189; green, 16; blue, 224 }  ,fill opacity=0.35 ] (266.1,81.06) -- (165.61,34.03) -- (164.98,178.29) -- cycle ;
\draw  [draw opacity=0][fill={rgb, 255:red, 189; green, 16; blue, 224 }  ,fill opacity=0.5 ] (130.16,107.05) -- (164.98,178.29) -- (164,34.12) -- cycle ;
\draw [color={rgb, 255:red, 0; green, 0; blue, 0 }  ,draw opacity=1 ][line width=1.3]    (507.9,106.41) -- (455.21,178.93) ;
\draw [color={rgb, 255:red, 0; green, 0; blue, 0 }  ,draw opacity=1 ][line width=1.3]    (507,106.48) -- (454.31,32.96) ;
\draw [line width=1.3]    (164.31,34.96) -- (164.31,180) ;
\draw [line width=1.3]  [dash pattern={on 5.63pt off 4.5pt}]  (129,106.73) -- (266,81.23) ;
\draw [color={rgb, 255:red, 0; green, 0; blue, 0 }  ,draw opacity=0.7 ] [dash pattern={on 4.5pt off 4.5pt}]  (79.06,152.3) -- (266,81.23) ;
\draw [color={rgb, 255:red, 0; green, 0; blue, 0 }  ,draw opacity=0.7 ] [dash pattern={on 4.5pt off 4.5pt}]  (79.06,62.66) -- (266,81.23) ;
\draw [color={rgb, 255:red, 0; green, 0; blue, 0 }  ,draw opacity=0.7 ][line width=0.75]    (217,107.48) -- (266,81.23) ;
\draw [line width=1.3]    (164.31,180) -- (266,81.23) ;
\draw [line width=1.3]    (164,34.12) -- (266,81.23) ;
\draw  [color={rgb, 255:red, 0; green, 0; blue, 0 }  ,draw opacity=0.7 ] (217.9,107.41) -- (165.21,179.93) -- (79.96,152.23) -- (79.96,62.59) -- (165.21,34.89) -- cycle ;
\draw [color={rgb, 255:red, 0; green, 0; blue, 0 }  ,draw opacity=0.7 ]   (79.06,62.66) -- (129,106.73) ;
\draw [line width=1.3]    (129,106.73) -- (164.31,180) ;
\draw [color={rgb, 255:red, 0; green, 0; blue, 0 }  ,draw opacity=1 ][line width=1.3]    (129,106.73) -- (164.31,34.96) ;
\draw [color={rgb, 255:red, 0; green, 0; blue, 0 }  ,draw opacity=0.7 ]   (129,106.73) -- (79.06,152.3) ;
\draw  [fill={rgb, 255:red, 9; green, 158; blue, 134 }  ,fill opacity=1 ] (125,106.73) .. controls (125,104.52) and (126.79,102.73) .. (129,102.73) .. controls (131.21,102.73) and (133,104.52) .. (133,106.73) .. controls (133,108.94) and (131.21,110.73) .. (129,110.73) .. controls (126.79,110.73) and (125,108.94) .. (125,106.73) -- cycle ;
\draw  [fill={rgb, 255:red, 9; green, 158; blue, 134 }  ,fill opacity=1 ] (160.31,34.96) .. controls (160.31,32.75) and (162.1,30.96) .. (164.31,30.96) .. controls (166.52,30.96) and (168.31,32.75) .. (168.31,34.96) .. controls (168.31,37.17) and (166.52,38.96) .. (164.31,38.96) .. controls (162.1,38.96) and (160.31,37.17) .. (160.31,34.96) -- cycle ;
\draw  [fill={rgb, 255:red, 9; green, 158; blue, 134 }  ,fill opacity=1 ] (213,107.48) .. controls (213,105.27) and (214.79,103.48) .. (217,103.48) .. controls (219.21,103.48) and (221,105.27) .. (221,107.48) .. controls (221,109.69) and (219.21,111.48) .. (217,111.48) .. controls (214.79,111.48) and (213,109.69) .. (213,107.48) -- cycle ;
\draw  [fill={rgb, 255:red, 9; green, 158; blue, 134 }  ,fill opacity=1 ] (160.31,180) .. controls (160.31,177.79) and (162.1,176) .. (164.31,176) .. controls (166.52,176) and (168.31,177.79) .. (168.31,180) .. controls (168.31,182.21) and (166.52,184) .. (164.31,184) .. controls (162.1,184) and (160.31,182.21) .. (160.31,180) -- cycle ;
\draw  [fill={rgb, 255:red, 9; green, 158; blue, 134 }  ,fill opacity=1 ] (75.06,152.3) .. controls (75.06,150.09) and (76.85,148.3) .. (79.06,148.3) .. controls (81.27,148.3) and (83.06,150.09) .. (83.06,152.3) .. controls (83.06,154.51) and (81.27,156.3) .. (79.06,156.3) .. controls (76.85,156.3) and (75.06,154.51) .. (75.06,152.3) -- cycle ;
\draw  [fill={rgb, 255:red, 9; green, 158; blue, 134 }  ,fill opacity=1 ] (75.06,62.66) .. controls (75.06,60.45) and (76.85,58.66) .. (79.06,58.66) .. controls (81.27,58.66) and (83.06,60.45) .. (83.06,62.66) .. controls (83.06,64.87) and (81.27,66.66) .. (79.06,66.66) .. controls (76.85,66.66) and (75.06,64.87) .. (75.06,62.66) -- cycle ;
\draw  [fill={rgb, 255:red, 9; green, 158; blue, 134 }  ,fill opacity=1 ] (262,81.23) .. controls (262,79.02) and (263.79,77.23) .. (266,77.23) .. controls (268.21,77.23) and (270,79.02) .. (270,81.23) .. controls (270,83.44) and (268.21,85.23) .. (266,85.23) .. controls (263.79,85.23) and (262,83.44) .. (262,81.23) -- cycle ;
\draw [line width=1.3]    (454.31,33.96) -- (454.31,179) ;
\draw [color={rgb, 255:red, 0; green, 0; blue, 0 }  ,draw opacity=0.7 ] [dash pattern={on 4.5pt off 4.5pt}]  (419,105.73) -- (556,80.23) ;
\draw [color={rgb, 255:red, 0; green, 0; blue, 0 }  ,draw opacity=0.7 ] [dash pattern={on 4.5pt off 4.5pt}]  (369.06,151.3) -- (556,80.23) ;
\draw [color={rgb, 255:red, 0; green, 0; blue, 0 }  ,draw opacity=0.7 ] [dash pattern={on 4.5pt off 4.5pt}]  (369.06,61.66) -- (556,80.23) ;
\draw [line width=1.3]    (507,106.48) -- (556,80.23) ;
\draw [line width=1.3]    (454.31,179) -- (556,80.23) ;
\draw [line width=1.3]    (454,33.12) -- (556,80.23) ;
\draw  [color={rgb, 255:red, 0; green, 0; blue, 0 }  ,draw opacity=0.7 ] (507.9,106.41) -- (455.21,178.93) -- (369.96,151.23) -- (369.96,61.59) -- (455.21,33.89) -- cycle ;
\draw [color={rgb, 255:red, 0; green, 0; blue, 0 }  ,draw opacity=0.7 ]   (369.06,61.66) -- (419,105.73) ;
\draw [color={rgb, 255:red, 0; green, 0; blue, 0 }  ,draw opacity=0.7 ][line width=0.75]    (419,105.73) -- (454.31,179) ;
\draw [color={rgb, 255:red, 0; green, 0; blue, 0 }  ,draw opacity=0.7 ][line width=0.75]    (419,105.73) -- (454.31,33.96) ;
\draw [color={rgb, 255:red, 0; green, 0; blue, 0 }  ,draw opacity=0.7 ]   (419,105.73) -- (369.06,151.3) ;
\draw  [fill={rgb, 255:red, 9; green, 158; blue, 134 }  ,fill opacity=1 ] (415,105.73) .. controls (415,103.52) and (416.79,101.73) .. (419,101.73) .. controls (421.21,101.73) and (423,103.52) .. (423,105.73) .. controls (423,107.94) and (421.21,109.73) .. (419,109.73) .. controls (416.79,109.73) and (415,107.94) .. (415,105.73) -- cycle ;
\draw  [fill={rgb, 255:red, 9; green, 158; blue, 134 }  ,fill opacity=1 ] (450.31,33.96) .. controls (450.31,31.75) and (452.1,29.96) .. (454.31,29.96) .. controls (456.52,29.96) and (458.31,31.75) .. (458.31,33.96) .. controls (458.31,36.17) and (456.52,37.96) .. (454.31,37.96) .. controls (452.1,37.96) and (450.31,36.17) .. (450.31,33.96) -- cycle ;
\draw  [fill={rgb, 255:red, 9; green, 158; blue, 134 }  ,fill opacity=1 ] (503,106.48) .. controls (503,104.27) and (504.79,102.48) .. (507,102.48) .. controls (509.21,102.48) and (511,104.27) .. (511,106.48) .. controls (511,108.69) and (509.21,110.48) .. (507,110.48) .. controls (504.79,110.48) and (503,108.69) .. (503,106.48) -- cycle ;
\draw  [fill={rgb, 255:red, 9; green, 158; blue, 134 }  ,fill opacity=1 ] (450.31,179) .. controls (450.31,176.79) and (452.1,175) .. (454.31,175) .. controls (456.52,175) and (458.31,176.79) .. (458.31,179) .. controls (458.31,181.21) and (456.52,183) .. (454.31,183) .. controls (452.1,183) and (450.31,181.21) .. (450.31,179) -- cycle ;
\draw  [fill={rgb, 255:red, 9; green, 158; blue, 134 }  ,fill opacity=1 ] (365.06,151.3) .. controls (365.06,149.09) and (366.85,147.3) .. (369.06,147.3) .. controls (371.27,147.3) and (373.06,149.09) .. (373.06,151.3) .. controls (373.06,153.51) and (371.27,155.3) .. (369.06,155.3) .. controls (366.85,155.3) and (365.06,153.51) .. (365.06,151.3) -- cycle ;
\draw  [fill={rgb, 255:red, 9; green, 158; blue, 134 }  ,fill opacity=1 ] (365.06,61.66) .. controls (365.06,59.45) and (366.85,57.66) .. (369.06,57.66) .. controls (371.27,57.66) and (373.06,59.45) .. (373.06,61.66) .. controls (373.06,63.87) and (371.27,65.66) .. (369.06,65.66) .. controls (366.85,65.66) and (365.06,63.87) .. (365.06,61.66) -- cycle ;
\draw  [fill={rgb, 255:red, 9; green, 158; blue, 134 }  ,fill opacity=1 ] (552,80.23) .. controls (552,78.02) and (553.79,76.23) .. (556,76.23) .. controls (558.21,76.23) and (560,78.02) .. (560,80.23) .. controls (560,82.44) and (558.21,84.23) .. (556,84.23) .. controls (553.79,84.23) and (552,82.44) .. (552,80.23) -- cycle ;

\draw (81.5,99.13) node [anchor=north west][inner sep=0.75pt]  [color={rgb, 255:red, 0; green, 0; blue, 0 }  ,opacity=0.7 ]  {$C_{\sigma } '$};
\draw (113.5,143.13) node [anchor=north west][inner sep=0.75pt]  [color={rgb, 255:red, 0; green, 0; blue, 0 }  ,opacity=0.7 ]  {$C_{\tau } '$};
\draw (109.5,66.13) node [anchor=north west][inner sep=0.75pt]  [color={rgb, 255:red, 0; green, 0; blue, 0 }  ,opacity=0.7 ]  {$C_{\rho } '$};
\draw (114,117.13) node [anchor=north west][inner sep=0.75pt]    {$r_{1} '$};
\draw (51,144.13) node [anchor=north west][inner sep=0.75pt]  [color={rgb, 255:red, 0; green, 0; blue, 0 }  ,opacity=0.7 ]  {$r_{2} '$};
\draw (50,54.13) node [anchor=north west][inner sep=0.75pt]  [color={rgb, 255:red, 0; green, 0; blue, 0 }  ,opacity=0.7 ]  {$r_{3} '$};
\draw (170,16.13) node [anchor=north west][inner sep=0.75pt]    {$r_{4} '$};
\draw (219,110.88) node [anchor=north west][inner sep=0.75pt]  [color={rgb, 255:red, 0; green, 0; blue, 0 }  ,opacity=0.7 ]  {$r_{5} '$};
\draw (166.31,183.4) node [anchor=north west][inner sep=0.75pt]    {$r_{6} '$};
\draw (276,72.13) node [anchor=north west][inner sep=0.75pt]    {$r_{7} '$};
\draw (142.5,80.13) node [anchor=north west][inner sep=0.75pt]    {$C_{1}'$};
\draw (372.5,98.13) node [anchor=north west][inner sep=0.75pt]  [color={rgb, 255:red, 0; green, 0; blue, 0 }  ,opacity=0.7 ]  {$C_{\sigma } '$};
\draw (403.5,142.13) node [anchor=north west][inner sep=0.75pt]  [color={rgb, 255:red, 0; green, 0; blue, 0 }  ,opacity=0.7 ]  {$C_{\tau } '$};
\draw (399.5,65.13) node [anchor=north west][inner sep=0.75pt]  [color={rgb, 255:red, 0; green, 0; blue, 0 }  ,opacity=0.7 ]  {$C_{\rho } '$};
\draw (406,114.13) node [anchor=north west][inner sep=0.75pt]  [color={rgb, 255:red, 0; green, 0; blue, 0 }  ,opacity=0.7 ]  {$r_{1} '$};
\draw (341,143.13) node [anchor=north west][inner sep=0.75pt]  [color={rgb, 255:red, 0; green, 0; blue, 0 }  ,opacity=0.7 ]  {$r_{2} '$};
\draw (340,53.13) node [anchor=north west][inner sep=0.75pt]  [color={rgb, 255:red, 0; green, 0; blue, 0 }  ,opacity=0.7 ]  {$r_{3} '$};
\draw (460,16.13) node [anchor=north west][inner sep=0.75pt]    {$r_{4} '$};
\draw (509,109.88) node [anchor=north west][inner sep=0.75pt]    {$r_{5} '$};
\draw (456.31,182.4) node [anchor=north west][inner sep=0.75pt]    {$r_{6} '$};
\draw (566,71.13) node [anchor=north west][inner sep=0.75pt]    {$r_{7} '$};
\draw (459.5,76.13) node [anchor=north west][inner sep=0.75pt]    {$C_{2}'$};

\end{tikzpicture}
    \end{center}

    \caption{Representation of the relevant cones of $\trop I'$, courtesy of F. Zaffalon.}
    \label{fig : conesNE}
\end{figure}
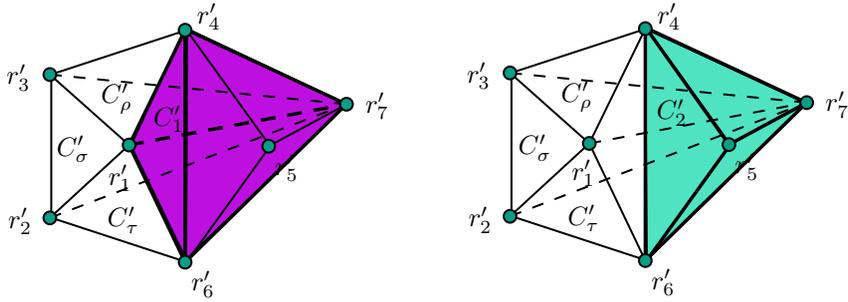

In \cref{fig : conesNE}, the relevant cones in the new embedding are depicted. This picture clarifies what is happening: we get a different subdivision of the bipyramid from \cref{fig : conesGr36}. The three non-prime cones making up the bipyramid in the original embedding have turned into two prime cones making up the same bipyramid. This also explains the existence of the prime cone $C_2'$, which projects onto all three of the non-prime cones, but is not adjacent to $C'_\rho$ or $C_\tau'$.

\subsection{Applying the wall crossing theorem}

We can now apply the wall crossing theorem to the cones we obtained. This will also resolve the issue with the choice of $C_1 $ and $C_2$, as we show below that the polytope $\Delta$ obtained in the original paper can be reconstructed from $C_1'$.

We only work out the wall-crossing for $C_{\rho}'$ and $C_1'$. First we have to choose 10 interior vectors of $C_{\rho}'\cap C'$ (we work in the full tropicalization here and do not quotient out the lineality space). Then we arrange them in a matrix $M_{C_{\rho}'\cap C_1'}$ as row vectors. Recall that we want the top row of the matrix to be identically one in order to simplify computations of the Newton-Okounkov body, which is possible as our ideal is homogeneous. In our computations we chose $M_{C_{\rho}'\cap C_1'}$ as follows: 

\[M_{C_{\rho}'\cap C_1'} = \left(\!\scalemath{0.6}{\begin{array}{cccccccccccccccccccccc}
1&1&1&1&1&1&1&1&1&1&1&1&1&1&1&1&1&1&1&1&1&1\\
0&0&0&-1&0&1&-1&0&0&-1&1&2&0&1&1&0&4&2&3&2&0&3\\
0&0&-1&0&0&0&0&-1&1&-1&1&1&1&0&2&0&3&3&3&2&0&3\\
0&-1&0&0&0&2&1&0&1&1&-1&1&0&-1&0&0&3&2&4&2&0&3\\
-1&0&0&0&0&1&0&1&2&1&-1&0&-1&0&1&0&3&2&3&3&0&3\\
1&1&1&1&1&2&1&1&2&1&0&1&0&0&1&0&3&2&3&2&1&3\\
1&1&1&1&1&2&1&1&2&1&0&1&0&0&1&0&3&2&3&2&0&4\\
0&0&0&0&0&2&0&0&1&0&0&2&0&0&1&0&4&2&4&2&0&4\\
0&0&0&0&0&1&0&0&1&0&0&1&0&0&1&0&4&3&4&3&0&4\\
0&0&0&0&0&1&0&0&2&0&0&1&0&0&2&0&4&3&4&3&0&4
\end{array}}\!\right) \]

Then we pick an interior vector of $C_{\rho}'$ and add it to the bottom of $M_{C_{\rho}'\cap C'}$ to obtain $M_\rho$ and analogously pick an interior vector of $C'$ to obtain $M$. We then obtain the corresponding Newton-Okounkov bodies $\Delta_\rho'$ and $\Delta'$ as the convex hull of the row vectors of $M_\rho$ and $M$ respectively, and the mutation can be obtained as described in section 3.4.3.

Note that the new-embedding adds two extra coordinates to the interior vectors relative to the original embedding, corresponding to the variables $Y$ and $h$. Therefore in order to return to the original embedding, one forgets the last two columns of the matrices $M_{C_{\rho}'\cap C'}$, $M_\rho$ or $M$. This means two vertices are deleted when returning to the original embedding.

However, the last column corresponds to the homogenization variable $h$, whose influence we would like to eliminate. This is because we want actually want to consider $I' \cap \langle h-1\rangle$ in order to obtain a variety isomorphic to $V(I_{3,6})$. This can be accomplished in the tropical side by intersecting $\trop(I')$ with the $h=0$ hyperplane (where we identify the last coordinate of the ambient space of $\trop(I')$ with $h$ by slight abuse of notation. Now our computations show that the vector $(0, ..., 0, 1, -1)$ is an element of the lineality space of $\trop(I')$. This means that in order to eliminate the $h$-column of our matrices, we can add the $h$-column to the $Y$ column and delete the $h$-column. The resulting added vertex is then $\frac{1}{2}$ times this sum of the $Y$- and $h$-columns, as we have to intersect with the $x_1 = 1$ hyperplane to get the Newton-Okounkov bodies.

Now by construction, the Newton-Okounkov body obtained by deleting the $Y$-vertex (and $h$-vertex) from $\Delta_\rho'$ is isomorphic to the matching field polytope $\Delta_\rho$. Deleting this vertex for $\Delta'$ does not make sense, as we could not obtain a Newton-Okounkov body through this construction from the original cone $C$ as it was non-prime. Instead, we compare $\Delta'$ with $\Delta$. More precisely, we compare $\Delta$ with $\Delta'$ where we removed the $h$-vertex and replaced it with the average of the $Y$- and $h$-vertices as discussed above. Remarkably, it is possible to computationally verify that these are isomorphic, which shows that our choice of $C_1$ and $C_2$ is correct.
%


%
  \graphicspath{{\chaptersdir/conclusion/image/}}%
\chapter{Conclusion \& further research}\label{ch:conclusion}


The goal of this thesis was to research a possible generalisation of the wall-crossing result for prime cones from \cite{WallCrossing} to the case of non-prime cones. This would allow for algorithmic computation of tropical varieties that have resisted brute forcing computations or computations making use of symmetry such as the higher dimensional tropical Grassmannians.

With this goal in mind we proposed an application of an algorithm from \cite{Algorithm}, which constructs a new embedding of the original variety and changes the tropicalization. The objective is to change the tropicalization in such a way that the non-prime cone becomes prime, and the adjacent cones stay adjacent. This would allow wall-crossing through the non-prime cone by lifting to the new embedding, applying the wall-crossing theorem for prime cones, and finally projecting back down.

Our computations on $\Gr(3,6)$ provide promising results, showing our application of the algorithm can provide a new tropicalization of the desired form. Furthermore, our computations neatly explain an unexpected sequence of mutations obtained in \cite{MatchingFields}, which served as the inspiration for this thesis. Furthermore we expect that our results will hold for the 14 other bipyramids of $\Gr(3,6)$.



It should be noted that while these results are promising, they are limited to $\Gr(3,6)$. In order to apply our approach more generally, more research is required into necessary and sufficient conditions for this algorithm to work.

Preservation of the adjacency of cones when computing the new embedding is one of the main concerns. This is not a problem we encounter in our computations for $\Gr(3,6)$, which is promising, but the same may not hold in other cases, for example if one needs to add multiple binomials, or apply the embedding algorithm multiple times. In general, studying the way this algorithm changes a tropicalization may prove worthwhile. For example, our computations show that the new tropicalization contains more maximal cones, and that multiple cones may project onto the same cone of the original embedding. 

Another concern is termination of the algorithm to compute the new embedding. While it has been shown that the algorithm does not always terminate in \cite{ilten2020khovanskiifinite}, explicit conditions for finishing remain unknown. For this purpose, it may be interesting to look at this algorithm from the more algebraic viewpoint of Khovanskii bases and valuations. 

Furthermore, this method still lacks the ability to calculate points of non-prime cones. While in the context of toric degenerations this is not a problem, in tropical geometry the whole tropical variety is of importance, not just the prime cones. Our results show that it may be possible to calculate the prime cones through the use of mutations, but non-prime cones remain unreachable. A more thorough understanding of how Algorithm 1 affects non-prime cones could function as a stepping stone for future research to reaching these non-prime cones.

Finally, it may yield to compute more examples and study the arising mutations to or through non-prime cones. Mutations between prime cones have known descriptions as more explicit operations on the underlying combinatorial descriptions of tropicalizations, such as for the space of phylogenetic trees for $\Gr(2,n)$. Studying our approach with these descriptions in mind could lead to insights into generating wall crossing maps without knowing the image polytope. This builds toward the end goal of implementing an algorithm allowing to compute an entire topicalization starting from a single known cone.

Overall our results show that non-prime wall crossing can correspond to mutations of Newton-Okounkov bodies, similar to wall-crossing for prime cones. Our results are promising, but more research is required in order to achieve the goal of implementing an algorithm to compute tropical varieties by wall-crossing.

\cleardoublepage


%




\backmatter

\includebibliography
\bibliographystyle{acm}
\bibliography{allpapers.bib}
\instructionsbibliography



\makebackcoverXII

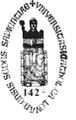
\end{document}